\documentclass[12pt,draft]{amsart}
\usepackage{amsmath,amsthm,latexsym,amscd,amsbsy,amssymb,pb-diagram}
 \relax

\chardef\bslash=`\\ 

\makeatletter
\def\verbatim{\interlinepenalty\@M \@verbatim
  \leftskip\@totalleftmargin\advance\leftskip2pc
  \frenchspacing\@vobeyspaces \@xverbatim}
\makeatother
\makeatletter
  \def\dgt@k{\dg@DX=-3 \dg@DY=2 \dg@SIZE=3} 
\makeatother

\makeatletter
  \def\dgt@kk{\dg@DX=3 \dg@DY=-1 \dg@SIZE=3}%
\makeatother

\theoremstyle{plain}
\newtheorem{thm}{Theorem}[section]
\newtheorem{cor}[thm]{Corollary}
\newtheorem{lem}[thm]{Lemma}
\newtheorem{pro}[thm]{Proposition}

\theoremstyle{definition}
\newtheorem{rem}[thm]{Remark}
\newtheorem{defin}[thm]{Definition}

\newcommand{\bb}{\bigstar_{\mathbb C}}

\numberwithin{equation}{section}

\begin{document}


\title[Uncountable direct systems and a characterization of
non-separable projective $C^{\ast}$-algebras]{Uncountable
direct systems and a characterization of non-separable
projective $C^{\ast}$-algebras}
\author{Alex Chigogidze}
\address{Department of Mathematics and Statistics,
University of Saskatche\-wan,
McLean Hall, 106 Wiggins Road, Saskatoon, SK, S7N 5E6,
Canada}
\email{chigogid@math.usask.ca}
\thanks{Author was partially supported by NSERC research grant.}



\begin{abstract}{We introduce the concept of a direct
$C_{\omega}^{\ast}$-system and show that every non-separable
unital $C^{\ast}$-algebra is the limit of essentially
unique direct $C_{\omega}^{\ast}$-system. This result is
then applied to the problem of characterization of projective
unital $C^{\ast}$-algebras. It is shown that a non-separable unital
$C^{\ast}$-algebra $X$ of density $\tau$ is projective if and
only if it is the limit of a well ordered direct system
${\mathcal S}_{X} = \{ X_{\alpha}, i_{\alpha}^{\alpha +1},
\alpha < \tau \}$ of length $\tau$, consisting of unital
projective $C^{\ast}$-subalgebras $X_{\alpha}$ of $X$ and doubly
projective homomorphisms (inclusions)
$i_{\alpha}^{\alpha +1} \colon X_{\alpha} \to X_{\alpha +1}$,
$\alpha < \tau$, so that $X_{0}$ is separable and each
$i_{\alpha}^{\alpha +1}$, $\alpha < \tau$, has a separable
type. In addition we show that a doubly projective homomorphism
$f \colon X \to Y$ of unital projective $C^{\ast}$-algebras has a
separable type
if and only if there exists a pushout diagram
\[
\begin{CD}
X @>f>> Y\\
@A{p}AA @AA{q}A\\
X_{0} @>f_{0}>> Y_{0},
\end{CD}
\]

\noindent where $X_{0}$ and $Y_{0}$ are separable unital projective
$C^{\ast}$-algebras and the homomorphisms $i_{0} \colon X _{0} \to Y_{0}$,
$p \colon X_{0} \to X$ and $q \colon Y_{0} \to Y$ are
doubly projective. These two results
provide a complete characterization of non-separable
projective unital $C^{\ast}$-algebras in terms of separable ones.}
\end{abstract}

\maketitle
\markboth{A.~Chigogidze}{Uncountable direct systems and a
characterization of non-separable projective $C^{\ast}$-algebras}

\section{Introduction}\label{S:intro}

The concept of the direct system in the $C^{\ast}$-algebra
theory has been successfully used in a wide range of situations.
While playing an important role in different constructions and
proofs of various statements, direct systems have been used as a
tool of introducing new concepts as well. But perhaps the most
significant demonstration of the power of direct systems is a
possibility of investigation of complicated $C^{\ast}$-algebras
by means of their approximation (via direct systems) by simpler
$C^{\ast}$-algebras. Such an approach is standard in almost any
category which possesses direct (or dually, inverse) systems.
Difficulties in systematic implementation of such a method have
variety of sources. If, for instance, we wish to investigate
a particular property of an arbitrarily given $C^{\ast}$-algebra $X$ by analyzing a randomly 
taken direct system
${\mathcal S}_{X}$, the limit of which is isomorphic to $X$, then we immediately face the 
fundamental problem of choice. Is the information encoded in the direct system 
${\mathcal S}_{X}$ relevant to the property of its limit under consideration? 
Does there exist a direct system with the same limit which is better designed for detecting 
that property? Are there effective ways of finding such a system? In other words, if two 
direct systems ${\mathcal S}_{X} = \{ X_{\alpha}, i_{\alpha}^{\beta}, A\}$ and 
${\mathcal S}_{Y} = \{ Y_{\alpha}, j_{\alpha}^{\beta}, A\}$ have isomorphic limits are 
these systems internally related to each other? Do they, for instance, contain isomorphic 
subsystems? Trivial examples show that the answer in general is negative. Two direct sequences
\[ {\mathcal S}_{0} = \left\{ C\left(\{ 0,1\}^{2n}\right) , C\left(\pi_{2n}^{2(n+1)}\right) , \omega\right\}\]

\noindent and
\[ {\mathcal S}_{1} = \left\{ C\left(\{ 0,1\}^{2n+1}\right) , C\left(\pi_{2n+1}^{2(n+1)+1}\right) ,\omega\right\} ,
\]

\noindent where $\{0,1\}$ is the two-point discrete space
and $\pi_{n}^{k} \colon \{ 0,1\}^{k} \to \{ 0,1\}^{n}$, $k > n$,
stands for the natural projection, obviously have the same
limit -- the $C^{\ast}$-algebra of continuous complex-valued
functions of the Cantor discontinuum -- but contain no isomorphic
subsequences whatsoever. 

In Section \ref{S:smooth} we introduce (Definition \ref{D:smooth})
the concept of direct $C_{\omega}^{\ast}$-system and prove
(Theorem \ref{T:spectral} and Proposition \ref{P:spectralis})
that if a unital $C^{\ast}$-algebra is represented
as the limit of two direct $C^{\ast}_{\omega}$-systems, then these systems
necessarily contain cofinal isomorphic subsystems. It is important to note
that every non-separable unital $C^{\ast}$-algebra is the limit of at least
one direct $C_{\omega}^{\ast}$-system (Proposition \ref{P:exists}). 
Therefore every non-separable unital $C^{\ast}$-algebra $X$ admits
essentially unique
direct $C_{\omega}^{\ast}$-system
${\mathcal S}_{X} = \{ X_{\alpha}, i_{\alpha}^{\beta}, A\}$ and
we conclude that any
information about $X$
is contained in ${\mathcal S}_{X}$. The remaining problem of recovering
such an information is, generally speaking, still quite challenging,
but has an explicit technical, and not a philosophical, nature.
An effective method of searching for such an information is based
on Proposition \ref{P:search}.

Actually Theorem \ref{T:spectral} states much more than it might
seem to be the case. Not only it states, as was indicated above,
that every two direct $C_{\omega}^{\ast}$-systems with isomorphic
limits contain isomorphic cofinal subsystems, but it essentially
guarantees that any homomorphism
$f \colon \varinjlim{\mathcal S}_{X} \to \varinjlim{\mathcal S}_{Y}$
between the limits of two direct $C_{\omega}^{\ast}$-systems
${\mathcal S}_{X} = \{ X_{\alpha}, i_{\alpha}^{\beta}, A\}$ and
${\mathcal S}_{Y} = \{ Y_{\alpha}, j_{\alpha}^{\beta}, A\}$ is
itself the limit $f = \varinjlim\{ f_{\alpha}; \alpha \in A_{f}\}$
of a certain morphism
$\{ f_{\alpha} \colon X_{\alpha} \to Y_{\alpha}, A_{f}\}
\colon {\mathcal S}_{X}|A_{f} \to {\mathcal S}_{Y}|A_{f}$,
consisting of ``level" homomorphisms,
between cofinal subsystems of the given ones. Such a phenomenom,
as was indicated above, is not possible for direct sequences. 

We apply the above outlined results to the problem of
characterization of non-separable projective unital
$C^{\ast}$-algebras in terms of separable ones. Here is the scheme we
follow. First we show (Lemma \ref{L:first}) that any non-separable projective
unital $C^{\ast}$-algebra $X$ is the limit of a direct
$C_{\omega}^{\ast}$-system
${\mathcal S}_{X} = \{ X_{\alpha}, i_{\alpha}^{\beta}, A\}$,
where $X_{\alpha}$'s, $\alpha \in A$, are separable projective
unital $C^{\ast}$-subalgebras of $X$ and the unital $\ast$-homomorphisms
$i_{\alpha}^{\beta} \colon X_{\alpha} \to X_{\beta}$, $\alpha \leq \beta$,
$\alpha ,\beta \in A$, are inclusions. It should be pointed out here
that the converse of this fact fails to be true, i.e. there does exist a
non-separable non-projective
unital $C^{\ast}$-algebra which is the limit of a direct
$C_{\omega}^{\ast}$-system consisting of separable and
projective unital $C^{\ast}$-subalgebras. This is how we arrive to
the necessity of analyzing inclusion homomorphisms $i_{\alpha}^{\beta}$.
What kind of property of these inclusion homomorphisms must be
present in order to guarantee that the limit of a direct system
${\mathcal S}_{X} = \{ X_{\alpha}, i_{\alpha}^{\beta}, A\}$,
consisting of projective $C^{\ast}$-subalgebras, is projective?
We especially emphasize this step because it is a crucial ingredient
of a typical argument based on Theorem
\ref{T:spectral}. In our particular situation explanation
is simple. The concept of a projective object has an explicit
categorical nature and seems logical to anticipate that the
required property of inclusion homomorphisms is closely related
to it. Consequently it makes sense to examine what does the
projectivity of a unital $\ast$-homomorphism, considered as
an object of the category $\operatorname{Mor}({\mathcal C}_{1}^{\ast})$
of unital $\ast$-homomorphisms of unital $C^{\ast}$-algebras,
mean. It turns out (Proposition \ref{P:morphisms}) that
projective objects of the category
$\operatorname{Mor}({\mathcal C}_{1}^{\ast})$ are precisely doubly
projective unital $\ast$-homomorphisms in the sense of \cite{lope1}.

In Section \ref{S:doubly} we establish certain properties of doubly
projective homomorphisms and present two characterizations of
non-separable projective unital $C^{\ast}$-algebras -- one
(condition (b) of Theorem \ref{T:charact}) in terms of direct
$C^{\ast}_{\omega}$-systems and the other (condition (c) of
Theorem \ref{T:charact}) in terms of well ordered continuous
direct systems. The latter states that a non-separable unital
$C^{\ast}$-algebra $X$ of density $\tau$ is projective if and
only if it is the limit of a well ordered direct system
${\mathcal S}_{X} = \{ X_{\alpha}, i_{\alpha}^{\alpha +1},
\alpha < \tau \}$ of length $\tau$, consisting of unital
projective $C^{\ast}$-subalgebras $X_{\alpha}$ of $X$ and doubly
projective homomorphisms (inclusions)
$i_{\alpha}^{\alpha +1} \colon X_{\alpha} \to X_{\alpha +1}$,
$\alpha < \tau$, so that $X_{0}$ is separable and each
$i_{\alpha}^{\alpha +1}$, $\alpha < \tau$, has a separable
type (Definition \ref{D:septype}).

Obviously this result can not be accepted as the one providing a
satisfactory reduction of the non-separable case to the separable
one. Of course, everything is fine if the density of $X$ is
$\omega_{1}$ -- in such a case all $X_{\alpha}$'s,
$\alpha < \omega_{1}$, (and not only the very first one,
i.e. $X_{0}$) are indeed separable. But if the density of $X$
is greater than $\omega_{1}$, then all $X_{\alpha}$'s, with
$\alpha \geq \omega_{1}$, are non-separable. 

In order to achieve our final goal and complete the reduction, we,
in Section \ref{S:diagrams},
analyze doubly projective homomorphisms of separable type between
(generally speaking, non-separable)
projective unital $C^{\ast}$-algebras. A characterization of such
homomorphisms, which is recorded in Theorem \ref{T:q}
(see also Corollary \ref{C:q}), states that a doubly projective homomorphism
$f \colon X \to Y$ of projective unital $C^{\ast}$-algebras has a separable type
if and only if there exists a pushout diagram
\[
\begin{CD}
X @>f>> Y\\
@A{p}AA @AA{q}A\\
X_{0} @>f_{0}>> Y_{0},
\end{CD}
\]

\noindent where $X_{0}$ and $Y_{0}$ are separable unital projective
$C^{\ast}$-algebras and the homomorphisms $i_{0} \colon X _{0} \to Y_{0}$,
$p \colon X_{0} \to X$ and $q \colon Y_{0} \to Y$ are doubly projective.

Theorems \ref{T:charact} and \ref{T:q} together complete the required reduction.

Proofs of above statements are based on some properties of unital free products. These properties are undoubtedly known to the experts in the field. For the readers convenience we discuss them in Section \ref{S:free}.

\bigskip
\bigskip


\section{Preliminaries}\label{S:pre}
All $C^{\ast}$-algebras below are assumed to be unital and all
$\ast$-homomorphisms between unital $C^{\ast}$-algebras are also
unital. The category formed by such $C^{\ast}$-algebras and
homomorphisms is denoted by ${\mathcal C}_{1}^{\ast}$. The density
$d(X)$ of a $C^{\ast}$-algebra $X$ is the minimal cardinality
of dense subspaces (in a purely topological sense) of $X$.
Thus $d(X) \leq \omega$ ($\omega$ denotes the first infinite
cardinal number) means that $X$ is separable. The unital
$C^{\ast}$-algebra,
consisting of only one element, is denoted by ${\mathbf 0}$.
${\mathbb C}$ denotes the $C^{\ast}$-algebra of complex numbers.

\subsection{Set-theoretical facts}\label{SS:set}
For the reader's convenience we begin by presenting necessary
set-theoretic facts. Their complete proofs can be found
in \cite{book}.

Let $A$  be a partially ordered {\em directed set} (i.e.
for every two elements  $\alpha ,\beta \in A$  there exists
an element  $\gamma \in A$  such that  $\gamma \geq \alpha$ 
and  $\gamma \geq \beta$). We say that a subset
$A_1 \subseteq A$ of $A$ {\em majorates} another subset
$A_2 \subseteq A$ of $A$ if for each element $\alpha_2 \in A_2$
there exists an element $\alpha_1 \in A_1$ such that
$\alpha_1 \geq \alpha_2$. A subset which majorates $A$
is called {\em cofinal} in $A$. A subset of  $A$  is said to
be a {\em chain} if every two elements of it are comparable.
The symbol $\sup B$ , where  $B \subseteq A$, denotes the
lower upper bound of $B$ (if such an element exists in $A$).
Let now $\tau$ be an infinite cardinal number. A subset $B$
of $A$  is said to be $\tau$-{\em closed} in $A$ if for each chain
$C \subseteq B$, with ${\mid}C{\mid} \leq \tau$, we have
$\sup C \in B$, whenever the element $\sup C$ exists in $A$.
Finally, a directed set $A$ is said to be $\tau$-{\em complete}
if for each chain $B$ of elements of $A$ with
${\mid}C{\mid} \leq \tau$, there exists an element
$\sup C$ in $A$. 

The standard example of a $\tau$-complete set can be obtained
as follows. For an arbitrary set $A$ let $\exp A$ denote, as usual,
the collection of all subsets of $A$. There is a natural partial
order on $\exp A$: $A_1 \geq A_2$ if and only if $A_1 \supseteq A_2$.
With this partial order $\exp A$ becomes a directed set.
If we consider only those subsets of the set $A$ which have
cardinality $\leq \tau$, then the corresponding subcollection
of $\exp A$, denoted by $\exp_{\tau}A$, serves as a basic
example of a $\tau$-complete set.

\begin{pro}\label{P:3.1.1}
Let  $\{ A_{t} : t \in T \}$ be a collection of $\tau$-closed and
cofinal subsets of a $\tau$-complete set $A$. If
$\mid T\mid \leq \tau$, then the intersection
$\cap \{ A_{t}: t \in T \}$ is also cofinal
(in particular, non-empty) and $\tau$-closed in $A$ .
\end{pro}

\begin{cor}\label{C:3.1.2}
For each subset $B$, with  $\mid B \mid \leq \tau$, of a
$\tau$-complete set $A$ there exists an element $\gamma \in A$
such that  $\gamma \geq \beta$  for each  $\beta \in B$ .
\end{cor}

\begin{pro}[Spectral Search]\label{P:search}
Let  $A$  be a $\tau$-complete set, 
$L \subseteq A^2$, and suppose the following three
conditions are satisfied:
\begin{description}
\item[Existence] For each $\alpha \in A$ there exists
$\beta \in A$  such that  $(\alpha ,\beta ) \in L$.
\item[Majorantness] If  $(\alpha ,\beta ) \in L$  and
$\gamma \geq \beta$, then  $(\alpha ,\gamma ) \in L$.
\item[$\tau$-closeness] Let $\{ \alpha_{t} : t \in T \}$
be a chain in $A$ with $|T| \leq \tau$. If
$(\alpha_{t}, \beta ) \in L$ for some
$\beta \in A$ and each $t \in T$, then
$(\alpha ,\beta ) \in L$ where $\alpha =
\sup \{\alpha_{t} \colon t \in T \}$.
\end{description}
   Then the set of all  $L$-{\em reflexive} elements of 
$A$ (an element $\alpha \in A$ is $L$-reflexive if
$(\alpha ,\alpha ) \in L$)  is cofinal and $\tau$-closed in $A$.
\end{pro}

Various applications of the above set-theoretical
statements are presented in \cite[Chapter 8]{book}. 


\bigskip
\bigskip

\section{Direct systems of unital $C^{\ast}$-algebras}\label{S:smooth}
Let us recall definitions of some of the concepts related to the notion
of a direct system.

\subsection{Morphisms of direct systems}\label{SS:morphisms}
A direct system
${\mathcal S} = \{ X_{\alpha}, i_{\alpha}^{\beta}, A\}$ of
unital $C^{\ast}$-algebras consists of a partially ordered directed
indexing set $A$,
unital $C^{\ast}$-algebras $X_{\alpha}$, $\alpha \in A$, and 
unital $\ast$-homomorphisms
$i_{\alpha}^{\beta} \colon X_{\alpha} \to X_{\beta}$,
defined for each
pair of indexes $\alpha ,\beta \in A$ with $\alpha \leq \beta$,
and satisfying
the condition $i_{\alpha}^{\gamma} =
i_{\beta}^{\gamma}\circ i_{\alpha}^{\beta}$ for
each triple of indexes
$\alpha ,\beta ,\gamma \in A$ with $\alpha \leq \beta \leq \gamma$.
The limit unital $C^{\ast}$-algebra of the above direct
system is denoted by
$\varinjlim{\mathcal S}$. For each $\alpha \in A$ there
exists a unital $\ast$-homomorphism
$i_{\alpha} \colon X_{\alpha} \to \varinjlim{\mathcal S}$
which will be called
the $\alpha$-th limit homomorphism of $\mathcal S$.

If $A^{\prime}$ is a
directed subset of the indexing set $A$, then the subsystem
$\{ X_{\alpha}, i_{\alpha}^{\beta}, A^{\prime}\}$ of
${\mathcal S}$ is denoted ${\mathcal S}|A^{\prime}$. 

Suppose that we are given two direct systems (with the
same indexing set)
${\mathcal S}_{X} = \{ X_{\alpha}, i_{\alpha}^{\beta}, A\}$ and
${\mathcal S}_{Y} = \{ Y_{\alpha}, j_{\alpha}^{\beta}, A \}$
consisting of unital $C^{\ast}$-algebras and unital $\ast$-ho\-mo\-mor\-phisms.
A {\em morphism} 
\[ \{ f_{\alpha} \colon \alpha \in A\} \colon
{\mathcal S}_{X} \to {\mathcal S}_{Y}\]

\noindent of the system ${\mathcal S}_X$ into
the system ${\mathcal S}_Y$ is a collection
$\{ f_{\alpha} \colon \alpha \in A\}$ of
unital $\ast$-ho\-mo\-mor\-phisms $f_{\alpha} \colon X_{\alpha}
\rightarrow Y_{\alpha}$, defined for all $\alpha \in A$,
such that\\
$$j^{\beta}_{\alpha}\circ f_{\alpha} = f_{\beta}\circ
i^{\beta}_{\alpha},$$
whenever $\alpha , \beta \in A$ and $\alpha \leq \beta$.
In other words, we require (in the above situation) the
commutativity of the following diagram

\[
\begin{CD}
X_{\beta} @>f_{\beta}>> Y_{\beta}\\
@A{i_{\alpha}^{\beta}}AA @AA{j^{\beta}_{\alpha}}A\\
X_{\alpha} @>f_{\alpha}>> Y_{\alpha}
\end{CD}
\]

\bigskip   

Any morphism $ \{ f_{\alpha} \colon \alpha \in A \} 
\colon {\mathcal S}_X \rightarrow {\mathcal S}_Y $ induces the
unital $\ast$-homomorphism, called {\em the
limit homomorphism of the morphism},\\
\[\varinjlim~\{ f_{\alpha} \colon \alpha \in A\}
\colon \varinjlim~{\mathcal S}_X \rightarrow
\varinjlim~{\mathcal S}_Y \]
\noindent such that
$\varinjlim~\{ f_{\alpha} \colon \alpha \in A\} \circ i_{\alpha} =
j_{\alpha}\circ f_{\alpha}$ for each $\alpha \in A$.
This obviously means
that all diagrams of the form

\[
\begin{CD}
\varinjlim{\mathcal S}_{X}
@>\varinjlim\{ f_{\alpha}\colon \alpha \in A\}>>
\varinjlim{\mathcal S}_{X}\\
@A{i_{\alpha}}AA @AA{j_{\alpha}}A\\
X_{\alpha} @>f_{\alpha}>> Y_{\alpha}
\end{CD}
\]
\noindent commute.

In particular, if for a direct system
${\mathcal S} = \{ X_{\alpha}, i_{\alpha}^{\beta}, A\}$
of unital $C^{\ast}$-algebras
and for a unital $C^{\ast}$-algebra $Y$, we are given
unital $\ast$-homomorphisms
$f_{\alpha} \colon X_{\alpha} \to Y$ so that
$f_{\alpha} = f_{\beta}\circ i_{\alpha}^{\beta}$ for each
$\alpha ,\beta \in A$ with $\alpha \leq \beta$, then there exists
the unique unital $\ast$-homomorphism
$\varinjlim\{ f_{\alpha} \colon \alpha \in A\} \colon
\varinjlim{\mathcal S} \to Y$ such that
$f_{\alpha} = \varinjlim\{ f_{\alpha} \colon \alpha \in A\}
\circ i_{\alpha}$ for each $\alpha \in A$. To see this apply
the above observation to the trivial direct system
${\mathcal S}_{Y} = \{ Y_{\alpha}, j_{\alpha}^{\beta},A\}$,
where $Y_{\alpha} = Y$ and $j_{\alpha}^{\beta} =
\operatorname{id}_{Y}$ for each $\alpha ,\beta \in A$
with $\alpha \leq \beta$.

In the cases when all
homomorphisms $i_{\alpha}^{\beta} \colon X_{\alpha} \to X_{\beta}$
and $i_{\alpha} \colon X_{\alpha} \to \varinjlim{\mathcal S}$
are inclusions
we will sometimes identify $X_{\alpha}$ with its image
$i_{\alpha}(X_{\alpha})$ in
$\varinjlim{\mathcal S}$ and denote the corresponding
direct system shortly
by ${\mathcal S} = \{ X_{\alpha}, A\}$.

A direct system ${\mathcal S}_{X} = \{ X_{\alpha}, i_{\alpha}^{\alpha +1}, \tau\}$, the indexing set of which is an infinite cardinal number $\tau$, is called well ordered. We say that such a direct system is continuous if for each limit ordinal number $\beta < \tau$ the homomorphism
\[ \varinjlim\{ i_{\alpha}^{\beta} ; \alpha < \beta\} \colon \varinjlim\{ X_{\alpha}, i_{\alpha}^{\alpha +1}, \beta\} \to X_{\beta}\]

\noindent is an isomorphism.


\subsection{Direct $C_{\tau}^{\ast}$-systems of
$C^{\ast}$-algebras}\label{SS:spectral}

The concept of the direct $C_{\tau}^{\ast}$-system, introduced in the 
following definition, will be used below.

\begin{defin}\label{D:smooth}
Let $\tau \geq \omega$
be a cardinal number. A direct system
${\mathcal S} = \{ X_{\alpha}, i_{\alpha}^{\beta}, A\}$ of
unital $C^{\ast}$-algebras
and unital $\ast$-ho\-mo\-morp\-hisms is called a {\em direct
$C_{\tau}^{\ast}$-system} if the
following conditions are satisfied:
\begin{itemize}
\item[(a)]
$A$ is a $\tau$-complete set.
\item[(b)]
Density of $X_{\alpha}$ is at most $\tau$
(i.e. $d(X_{\alpha}) \leq \tau$), $\alpha \in A$.
\item[(c)]
The $\alpha$-th limit homomorphism
$i_{\alpha} \colon X_{\alpha} \to \varinjlim{\mathcal S}$ is
an injective $\ast$-ho\-mo\-mor\-phism for each $\alpha \in A$.
\item[(d)]
If $B = \{ \alpha_{t} \colon t \in T\} $ is a chain of elements of $A$ with
$|T| \leq \tau$ and $\alpha = \sup B$, then the limit
homomorphism
$\varinjlim\{ i_{\alpha_{t}}^{\alpha} \colon t \in T\}
\colon \varinjlim\left({\mathcal S}_{X}|B\right)
\to X_{\alpha}$ is an isomorphism.
\end{itemize}
\end{defin}

\begin{pro}\label{P:exists}
Let $\tau$ be an infinite cardinal number. Every unital
$C^{\ast}$-algebra $X$
can be represented as the limit of a direct
$C_{\tau}^{\ast}$-system
${\mathcal S}_{X} = \{ X_{\alpha}, i_{\alpha}^{\beta},
\exp_{\tau}d(X) \}$.
\end{pro}
\begin{proof}
If $d(X) \leq \tau$, then consider the direct $C_{\tau}^{\ast}$-system
${\mathcal S}_{X} = \{ X_{\alpha}, i_{\alpha}^{\beta},
\exp_{\tau}d(X) \}$, where $X_{\alpha} = X$ for each
$\alpha \in \exp_{\tau}d(X)$ and
$i_{\alpha}^{\beta} = \operatorname{id}_{X}$ for each
$\alpha ,\beta \in \exp_{\tau}d(X)$ with
$\alpha \leq \beta$.

If $d(X) > \tau$, then consider any subset $Y$ of $X$ such that
$\operatorname{cl}_{X}Y = X$ and $|Y| = d(X)$. Without loss
of generality we may assume that $Y$ contains the unit of $X$. Each
$\alpha \in \exp_{\tau}d(X)$ can obviously be identified with a subset
(denoted by the same letter $\alpha$) of $Y$ of cardinality $\leq \tau$.
Let $X_{\alpha}$ be the smallest $C^{\ast}$-subalgebra of $X$
containing $\alpha$.
If $\alpha ,\beta \in \exp_{\tau}d(X)$ and $\alpha \leq \beta$, then
$\alpha \subseteq \beta$ (as subsets of $Y$) and consequently
$X_{\alpha} \subseteq X_{\beta}$. This inclusion map is denoted
by $i_{\alpha}^{\beta} \colon X_{\alpha} \to X_{\beta}$. It
is easy to verify
that the collection
${\mathcal S}_{X} = \{ X_{\alpha}, i_{\alpha}^{\beta},
\exp_{\tau}d(X)\}$ is indeed
a direct $C_{\tau}^{\ast}$-system such that
$\varinjlim{\mathcal S}_{X} = X$.
\end{proof}

\begin{lem}\label{L:strong}
If ${\mathcal S}_{X} = \{ X_{\alpha}, i_{\alpha}^{\beta}, A\}$
is a direct $C_{\tau}^{\ast}$-system, then 
\[ \varinjlim{\mathcal S}_{X} = \bigcup\{
i_{\alpha}(X_{\alpha}) \colon \alpha \in A\} .\]
\end{lem}
\begin{proof}
Clearly $\displaystyle \bigcup\{ i_{\alpha}(X_{\alpha})
\colon \alpha \in A\}$
is dense in $\varinjlim{\mathcal S}_{X}$ (this fact
remains true for arbitrary direct
systems of $C^{\ast}$-algebras). Consequently, for any point
$x \in \varinjlim{\mathcal S}_{X}$ there exists a sequence
$\{ x_{n} \colon n \in \omega\}$, consisting of elements from
$\displaystyle \bigcup\{ i_{\alpha}(X_{\alpha})
\colon \alpha \in A\}$, such that $x = \lim\{ x_{n}
\colon n \in \omega\}$. For each $n \in \omega$ choose an index
$\alpha_{n} \in A$ such that
$x_{n} \in i_{\alpha_{n}}\left( X_{\alpha_{n}}\right)$.
By Corollary \ref{C:3.1.2}, there exists an index $\alpha \in A$
such
that $\alpha \geq \alpha_{n}$ for each $n \in \omega$.
Since $i_{\alpha_{n}} = i_{\alpha}\circ i_{\alpha_{n}}^{\alpha}$,
it follows that
\[ x_{n} \in i_{\alpha_{n}}\left( X_{\alpha_{n}}\right) =
i_{\alpha}\left( i_{\alpha_{n}}^{\alpha}\left(
X_{\alpha_{n}}\right)\right) \subseteq i_{\alpha}\left(
X_{\alpha}\right) \;\;\text{for each}\;\; n \in \omega .\]
Finally, since $i_{\alpha}\left( X_{\alpha}\right)$ is closed
in $\varinjlim{\mathcal S}_{X}$, it follows that 
\[ x = \lim\{ x_{n} \colon n \in \omega\} \in
i_{\alpha}\left( X_{\alpha}\right) .\]
\end{proof}

\begin{lem}\label{L:factor}
Let ${\mathcal S}_{X} = \{ X_{\alpha}, i_{\alpha}^{\beta}, A\}$
be a direct $C_{\tau}^{\ast}$-system and
$f \colon Y \to \varinjlim{\mathcal S}_{X}$
be a unital $\ast$-ho\-mo\-morp\-hism of a unital $C^{\ast}$-algebra $Y$
into the direct limit of ${\mathcal S}_{X}$. If $d(Y) \leq \tau$,
then there exist an index $\alpha \in A$ and a unital $\ast$-ho\-mo\-morp\-hism
$f_{\alpha} \colon Y \to X_{\alpha}$ such that
$f = i_{\alpha}\circ f_{\alpha}$. 
\end{lem}
\begin{proof}
Since $d(Y) \leq \tau$, there exists a dense subset
$Z = \{ z_{t} \colon t \in T\}$ of $Y$ such that
$|T| \leq \tau$. For each $t \in T$ there exists, by Lemma
\ref{L:strong},
an index $\alpha_{t} \in A$ such that
$f(z_{t}) \in i_{\alpha_{t}}\left(X_{\alpha_{t}}\right)$.
Since $A$ is a
$\tau$-complete (condition (a) of Definition \ref{D:smooth}),
there exists,
by Corollary \ref{C:3.1.2}, an index $\alpha \in A$ such that
$\alpha \geq \alpha_{t}$ for each $t \in T$. As in the proof of
Lemma \ref{L:strong} we can conclude that 
\[ f(Z) = f\left(\{ z_{t} \colon t \in T\}\right) =
\{ f(z_{t}) \colon t \in T\} \subseteq
i_{\alpha}\left( X_{\alpha}\right) .\]
Since $Z$ is dense in $Y$ and since
$i_{\alpha}\left( X_{\alpha}\right)$
is closed in $\varinjlim{\mathcal S}_{X}$ it follows that
\[ f(Y) = f\left(\operatorname{cl}_{Y}Z\right) \subseteq
\operatorname{cl}_{\varinjlim{\mathcal S_{X}}}f(Z) \subseteq
\operatorname{cl}_{\varinjlim{\mathcal S_{X}}}i_{\alpha}
\left( X_{\alpha}\right) =
i_{\alpha}\left( X_{\alpha}\right) .\]
By condition (c) of Definition \ref{D:smooth}, the $\alpha$-th limit
homomorphism $i_{\alpha}$ of the direct
$C_{\tau}^{\ast}$-system ${\mathcal S}_{X}$
is an injective unital $\ast$-ho\-mo\-morp\-hism. Thus the composition
$f_{\alpha} = i_{\alpha}^{-1}\circ f \colon Y \to X_{\alpha}$ is
a well defined unital
$\ast$-homomorphism. It only remains to note that
$i_{\alpha}\circ f_{\alpha} =
i_{\alpha}\circ i_{\alpha}^{-1}\circ f = f$, as required.
\end{proof}

The following statement is one of our main results.

\begin{thm}\label{T:spectral}
Let ${\mathcal S}_{X} = \{ X_{\alpha}, i_{\alpha}^{\beta}, A\}$
and ${\mathcal S}_{Y} = \{ Y_{\alpha}, j_{\alpha}^{\beta}, A\}$ be two
direct $C_{\tau}^{\ast}$-systems with the same indexing set $A$.
If $f \colon \varinjlim{\mathcal S}_{X} \to \varinjlim{\mathcal S}_{Y}$
is a unital $\ast$-ho\-mo\-mor\-phism between the limit $C^{\ast}$-algebras
of ${\mathcal S}_{X}$ and 
${\mathcal S}_{Y}$, then there exist a cofinal and $\tau$-closed subset
$A_{f} \subseteq A$ and a morphism 
\[ \{ f_{\alpha} \colon X_{\alpha} \to Y_{\alpha},
\alpha \in A_{f}\} \colon {\mathcal S}_{X}|A_{f} \to
{\mathcal S}_{Y}|A_{f} \]
such that $f = \varinjlim\{ f_{\alpha} \colon \alpha \in A_{f}\}$. 
\end{thm}
\begin{proof}
We perform the spectral search (see Proposition \ref{P:search})
with respect
to the relation $L_{f} \subseteq A^{2}$ which is defined as follows.
An ordered pair $(\alpha , \beta )$ of indeces is an element of
$L_{f}$ if and only if $\alpha \leq \beta$ and there
exists a unital $\ast$-ho\-mo\-morp\-hism
$f_{\alpha}^{\beta} \colon X_{\alpha} \to Y_{\beta}$ such that
$f\circ i_{\alpha} = j_{\beta}\circ f_{\alpha}^{\beta}$, i.e. if the
diagram

\[
\begin{CD}
\varinjlim{\mathcal S}_{X} @>f>> \varinjlim{\mathcal S}_{Y}\\
@A{i_{\alpha}}AA @AA{j_{\beta}}A\\
X_{\alpha} @>f_{\alpha}^{\beta}>> Y_{\beta}
\end{CD}
\]
\noindent commutes. Let us verify conditions of
Proposition \ref{P:search}.

{\em Existence}. For each $\alpha \in A$ we need to
find an index $\beta \in A$ such that
$(\alpha ,\beta ) \in L_{f}$. Indeed, according to condition
(b) of Definition \ref{D:smooth}, $d\left( X_{\alpha}\right) \leq \tau$. 
Consider the unital $\ast$-homomorphism $f\circ i_{\alpha}
\colon X_{\alpha} \to \varinjlim{\mathcal S}_{Y}$. By Lemma
\ref{L:factor}, there exist an index
$\beta \in A$ (which, without loss of generality, may
be assumed to be greater than $\alpha$) and a unital $\ast$-ho\-mo\-morp\-hism
$f_{\alpha}^{\beta} \colon X_{\alpha} \to Y_{\beta}$ such that
$f\circ i_{\alpha} = j_{\beta}\circ f_{\alpha}^{\beta}$. This
obviously means that $(\alpha ,\beta ) \in L_{f}$.

{\em Majorantness.} Let $(\alpha ,\beta ) \in L_{f}$ and
$\gamma \geq \beta$. In order to show that
$(\alpha ,\gamma ) \in L_{f}$,
consider the composition
$f_{\alpha}^{\gamma} = j_{\beta}^{\gamma} \circ
f_{\alpha}^{\beta} \colon X_{\alpha} \to Y_{\gamma}$, where the
unital $\ast$-ho\-mo\-morp\-hism $f_{\alpha}^{\beta} \colon X_{\alpha}\to Y_{\beta}$
is supplied by the condition $(\alpha ,\beta ) \in L_{f}$. Clearly
$j_{\gamma}\circ f_{\alpha}^{\gamma} =
j_{\gamma}\circ j_{\beta}^{\gamma}\circ
f_{\alpha}^{\beta} = j_{\alpha}\circ
f_{\alpha}^{\beta} = f\circ i_{\alpha}$.
This shows that $(\alpha ,\gamma ) \in L_{f}$.

{\em $\tau$-closeness}. Let $B = \{ \alpha_{t} \colon t \in T\}$
be a chain of elements in $A$ with $|T| \leq \tau$. Suppose that
$(\alpha_{t}, \beta ) \in L_{f}$ for some
$\beta \in A$ and each $t \in T$. We need to show
that $(\alpha ,\beta ) \in L_{f}$, where
$\alpha = \sup\{ \alpha_{t} \colon t \in T\}$.
First observe that if $\alpha_{t} \leq \alpha_{s}$ for $t,s \in T$, then
\[ j_{\beta}\circ f_{\alpha_{t}}^{\beta} =
f\circ i_{\alpha_{t}} =
f\circ i_{\alpha_{s}}\circ i_{\alpha_{t}}^{\alpha_{s}}
= j_{\beta}\circ f_{\alpha_{s}}^{\beta}\circ i_{\alpha_{t}}^{\alpha_{s}}.\]
Since, by condition (c) of Definition \ref{D:smooth}, the $\beta$-th
limit $\ast$-homomorphism $j_{\beta}$ of the direct system
${\mathcal S}_{Y}$ is injective, it follows that
$f_{\alpha_{t}}^{\beta} =
f_{\alpha_{s}}^{\beta}\circ i_{\alpha_{t}}^{\alpha_{s}}$. This means
that the collection $\{ f_{\alpha_{t}}^{\beta} \colon t \in T\}$
forms a morphism
of the subsystem ${\mathcal S}_{X}|B$ of the direct
$C_{\tau}^{\ast}$-system ${\mathcal S}_{X}$ into the
$C^{\ast}$-algebra $Y_{\beta}$. Consider (see
Subsection \ref{SS:morphisms}) the unital $\ast$-ho\-mo\-morp\-hism
\[ \varinjlim\{  f_{\alpha_{t}}^{\beta} \colon t \in T\}
\colon \varinjlim\left( {\mathcal S}_{X}|B\right) \to Y_{\beta} .\]
Finally, applying condition (d) of Definition \ref{D:smooth},
we define the unital $\ast$-ho\-mo\-morp\-hism
$f_{\alpha}^{\beta} \colon X_{\alpha} \to Y_{\beta}$
as the composition
\[ X_{\alpha} \xrightarrow{\left(\varinjlim\{
i_{\alpha_{t}}^{\alpha} \colon t \in T\}\right)^{-1}}
\varinjlim\left( {\mathcal S}_{X}|B\right)
\xrightarrow{\varinjlim\{  f_{\alpha_{t}}^{\beta}
\colon t \in T\}}Y_{\beta} .\] 
The straightforward verification shows that $f_{\alpha}^{\beta}$
indeed satisfies the required equality
$f\circ i_{\alpha} = j_{\beta}\circ f_{\alpha}^{\beta}$
and, consequently, $(\alpha ,\beta ) \in L_{f}$.

Now, by applying Proposition \ref{P:search}, we conclude that the set
$A_{f}$ of
$L_{f}$-reflexive elements is cofinal and $\tau$-closed in $A$.
Observe that an element $\alpha \in A$ is $L_{f}$-reflexive
if and only if there exists a unital $\ast$-ho\-mo\-morp\-hism
$f_{\alpha} \colon X_{\alpha} \to Y_{\alpha}$ such that
$f\circ i_{\alpha} = j_{\alpha}\circ f_{\alpha}$.

It follows from the above construction that the collection
\[ \{ f_{\alpha} \colon X_{\alpha} \to Y_{\alpha},
\alpha \in A_{f}\} \colon {\mathcal S}_{X}|A_{f} \to
{\mathcal S}_{X}|A_{f} \]
is indeed a morphism between the systems ${\mathcal S}_{X}|A_{f}$
and ${\mathcal S}_{Y}|A_{f}$
such that $f = \varinjlim\{ f_{\alpha} \colon \alpha \in A_{f}\}$. 
\end{proof}

\begin{pro}\label{P:spectralis}
If $f \colon \varinjlim{\mathcal S}_{X} \to
\varinjlim{\mathcal S}_{Y}$ is a unital$\ast$-isomorphism between the
limit $C^{\ast}$-algebras of direct $C_{\tau}^{\ast}$-systems
${\mathcal S}_{X} = \{ X_{\alpha}, i_{\alpha}^{\beta}, A\}$ and
${\mathcal S}_{Y} = \{ Y_{\alpha}, j_{\alpha}^{\beta}, A\}$ with
the same indexing set, then there exist a cofinal and $\tau$-closed
subset $A_{f} \subseteq A$ and a morphism 
\[ \{ f_{\alpha} \colon X_{\alpha} \to Y_{\alpha},
\alpha \in A_{f}\} \colon {\mathcal S}_{X}|A_{f} \to
{\mathcal S}_{Y}|A_{f}\]
such that $f = \varinjlim\{ f_{\alpha}\colon \alpha
\in A_{f}\}$ and $f_{\alpha}$ is a unital $\ast$-isomorphism
for each $\alpha \in A_{f}$.
\end{pro}
\begin{proof}
By Theorem \ref{T:spectral}, applied to the unital $\ast$-homomorphism
$f \colon \varinjlim{\mathcal S}_{X} \to
\varinjlim{\mathcal S}_{Y}$,
there exist a cofinal and $\tau$-closed subset
$\widetilde{A}_{f} \subseteq A$ and a morphism 
\[ \{ f_{\alpha} \colon X_{\alpha} \to Y_{\alpha},
\alpha \in \widetilde{A}_{f}\} \colon
{\mathcal S}_{X}|\widetilde{A}_{f} \to
{\mathcal S}_{Y}|\widetilde{A}_{f}\]
such that $f = \varinjlim\{ f_{\alpha}\colon \alpha
\in \widetilde{A}_{f}\}$.

Similarly, by Theorem \ref{T:spectral}, applied to
the unital $\ast$-homomorphism\\
$f^{-1} \colon \varinjlim{\mathcal S}_{Y} \to
\varinjlim{\mathcal S}_{X}$ (recall that $f$ is a unital $\ast$-isomorphism),
there exist a cofinal and $\tau$-closed subset
$\widetilde{A}_{f^{-1}} \subseteq A$ and a morphism 
\[ \{ g_{\alpha} \colon Y_{\alpha} \to X_{\alpha},
\alpha \in \widetilde{A}_{f^{-1}}\} \colon
{\mathcal S}_{Y}|\widetilde{A}_{f^{-1}} \to
{\mathcal S}_{X}|\widetilde{A}_{f^{-1}}\]
such that $f^{-1} = \varinjlim\{ g_{\alpha}\colon \alpha
\in \widetilde{A}_{f^{-1}}\}$.

By Proposition \ref{P:3.1.1}, the intersection
$A_{f} = \widetilde{A}_{f}\cap \widetilde{A}_{f^{-1}}$ is
still cofinal and $\tau$-closed subset of $A$. Note that for each
$\alpha \in A_{f}$ we have two unital $\ast$-homomorphisms
$f_{\alpha} \colon X_{\alpha} \to Y_{\alpha}$ and
$g_{\alpha} \colon Y_{\alpha} \to X_{\alpha}$ satisfying the equalities
$f\circ i_{\alpha} = j_{\alpha}\circ f_{\alpha}$ and
$f^{-1}\circ j_{\alpha} = i_{\alpha}\circ g_{\alpha}$. Consequently,
having also in mind condition (c) of Definition \ref{D:smooth}, we have
\begin{multline*}
 g_{\alpha}\circ f_{\alpha} =
g_{\alpha}\circ j_{\alpha}^{-1}\circ j_{\alpha} \circ f_{\alpha} =
g_{\alpha}\circ j_{\alpha}^{-1}\circ f \circ i_{\alpha} =
i_{\alpha}^{-1}\circ i_{\alpha} \circ g_{\alpha} =
j_{\alpha}^{-1} \circ f \circ i_{\alpha} =\\
 i_{\alpha}^{-1}\circ f^{-1}\circ j_{\alpha}\circ
j_{\alpha}^{-1}\circ f \circ i_{\alpha} =
i_{\alpha}^{-1}\circ f^{-1}\circ f \circ i_{\alpha} =
i_{\alpha}^{-1}\circ i_{\alpha} = \operatorname{id}_{X_{\alpha}}.
\end{multline*} 
Similarly, $f_{\alpha}\circ g_{\alpha} = \operatorname{id}_{Y_{\alpha}}$. 
This obviously means that both $f_{\alpha}$ and $g_{\alpha}$
are $\ast$-isomorphisms (inverses of each other).
\end{proof}


\bigskip
\bigskip

\section{Unital free products of unital
$C^{\ast}$-algebras and their direct
$C_{\omega}^{\ast}$-systems}\label{S:free}
Definition and various properties
of (amalgamated) free products
of $C^{\ast}$-algebras can be found in  \cite{brown},
\cite{cuntz}, \cite{effros}, \cite{blackadar1},
\cite{blackadar2}, \cite{pedersen2}. Below we consider only
the unital free products of unital
$C^{\ast}$-algebras (see \cite[Example 1.3(f)]{blackadar2}).

Recall that the unital free product of a collection
$\{ X_{t} \colon t \in T\}$ of unital $C^{\ast}$-algebras
(i.e. the amalgamated free product over the
common unit of $X_{t}$'s) is the 
unital $C^{\ast}$-algebra $\bb \{ X_{t} \colon t \in T\}$, together
with unital injective $\ast$-ho\-mo\-mor\-phisms
$\pi_{t} \colon X_{t} \to \bb \{ X_{t} \colon t \in T\}$,
$t \in T$, satisfying the following universality property:

\begin{itemize}
\item[$\bigl(\bigstar_{\mathbb C}\bigr)$]
for any unital $\ast$-homomorphisms
$f_{t} \colon X_{t} \to Y$, $t \in T$, into any
unital $C^{\ast}$-algebra $Y$, there exists a unique
unital $\ast$-homomorphism\\
$\bb \{ f_{t} \colon t \in T\} \colon \bb
\{X_{t} \colon t \in T\} \to Y$ such that
$f \circ \pi_{t} = f_{t}$, $t \in T$.
\end{itemize}

For the case $|T| = 2$ this universality property
of the unital free products is explicitly
stated by several authors (see, for instance, \cite[p. 81]{cuntz}
\cite[p. 156]{effros}, \cite[p. 89]{lope1},
\cite[2.2. Pushouts]{pedersen1}). The existence of unital free products
of infinite (uncountable)
collections of unital $C^{\ast}$-algebras is proved in 
\cite[Theorem 3.1]{blackadar1} (see also
\cite[Example 1.3(f)]{blackadar2}). Even though the property
$\bigl(\bigstar_{\mathbb C}\bigr)$ is not explicitly stated in \cite{blackadar1},
it can be extracted from the proofs provided there.

An alternative
approach for
establishing the property $\bigl(\bigstar_{\mathbb C}\bigr)$
for arbitrary $T$ is standard, although
less constructive. It is based on the following observation. The
unital direct product $\left( X_{1}\bb X_{2}, \pi_{1}, \pi_{2}\right)$
of two unital $C^{\ast}$-algebras $X_{1}$ and $X_{2}$ is precisely the
coproduct of the objects $X_{1}$ and $X_{2}$ in the category
${\mathcal C}_{1}$ (see \cite[p. 63]{maclane}). Also note that the category
${\mathcal C}_{1}$ is the category with the initial object -- namely,
the $C^{\ast}$-algebra ${\mathbb C}$. These two facts suffice
\cite[Section III.5]{maclane} to conclude the existence of unital
free products of finite collections of unital $C^{\ast}$-algebras.

For an arbitrary indexing set $T$, consider the directed set
$\exp_{<\omega}T$ of all finite subsets of $T$
with the natural partial order generated by the inclusion relation.
Next consider the direct system

\[ {\mathcal S} = \Bigl\{ \bb \{ X_{t} \colon t \in A\} ,
\pi_{A}^{B}, A, B \in \exp_{<\omega}T\Bigr\} ,\]
consisting of the unital free products
$\bb \{ X_{t} \colon t \in A\}$ of finite subcollections
and associated injective unital $\ast$-homomorphisms 

\[ \pi_{A}^{B} \colon \bb \{ X_{t} \colon t \in A\} \to
\bb \{ X_{t} \colon t \in B\} ,
A \subseteq B , A, B  \in \exp_{<\omega}T .\]

\noindent The limit $\varinjlim{\mathcal S}$ of this direct
system is in fact the unital free product of
the given collection. Unital injective
$\ast$-homomorphisms $\pi_{t}$, participating in
the definition of unital free products, are precisely
the $t$-th limit $\ast$-homomorphisms of the direct
system $\mathcal S$ (see Subsection \ref{SS:morphisms}).

To see that the property $\bigl(\bigstar_{\mathbb C}\bigr)$
is satisfied for so defined unital free products,
consider unital $\ast$-homomorphisms
\[ \bb \{ f_{t} \colon t \in A\} \colon
\bb \{ X_{t} \colon t \in A\} \to Y ,\]
uniquely defined (in the above discussed case of
the unital free products of finite collections of
unital $C^{\ast}$-algebras) for each finite subset
$A$ of $T$. It follows that 
\[  \bb \{ f_{t} \colon t \in A\} = 
\bb \{ f_{t} \colon t \in B\} \circ \pi_{A}^{B} \]
whenever $A, B \in \exp_{<\omega}T$ and $A \subseteq B$.
This guarantees
(see Subsection \ref{SS:morphisms}) that the limit
$\ast$-homomorphism 
\[ \bb\{ f_{t} \colon t \in T\} =
\varinjlim\{ \bb \{ f_{t} \colon t \in A\}
\colon A \in \exp_{<\omega}T\} \colon
\varinjlim{\mathcal S} \to Y\]
satisfies the required equalities
\[  \bb\{ f_{t} \colon t \in T\} \circ \pi_{t} =
f_{t} , t \in T .\]

We now state some of the properties of the unital free products
which will be needed in later sections.

\begin{lem}\label{L:2}
If $S \subseteq T$, then
$\displaystyle \bigstar_{{\mathbb C}}\{ X_{t}
\colon t \in T\}$ is canonically
isomorphic to\\ 
$\displaystyle \left(\bigstar_{{\mathbb C}}\{ X_{t}
\colon t \in S\}\right)
\bigstar_{{\mathbb C}}\left(\bigstar_{{\mathbb C}}\{
X_{t} \colon t \in T-S\}\right)$.
\end{lem}
\begin{proof}
Let 
\[\pi_{1} \colon \bigstar_{{\mathbb C}}\{ X_{t}
\colon t \in S\} \to \left(\bigstar_{{\mathbb C}}\{ X_{t}
\colon t \in S\}\right)
\bigstar_{{\mathbb C}}\left(\bigstar_{{\mathbb C}}\{
X_{t} \colon t \in T-S\}\right)\]

\noindent and

\[\pi_{2} \colon \bigstar_{{\mathbb C}}\{ X_{t}
\colon t \in T-S\} \to \left(\bigstar_{{\mathbb C}}\{ X_{t}
\colon t \in S\}\right)
\bigstar_{{\mathbb C}}\left(\bigstar_{{\mathbb C}}\{
X_{t} \colon t \in T-S\}\right)\]

\noindent denote the canonical inclusions
(see \cite[Theorem 3.1]{blackadar1}).

Let  
\[  \pi_{t}^{S} \colon X_{t} \to  \bigstar_{{\mathbb C}}\{ X_{t}
\colon t \in S\}\;\; ,\;
\pi_{t}^{T-S} \colon X_{t} \to  \bigstar_{{\mathbb C}}\{ X_{t}
\colon t \in T-S\}\]
\noindent and
\[  \pi_{t}^{T} \colon X_{t} \to  \bigstar_{{\mathbb C}}\{ X_{t}
\colon t \in T\}\]
\noindent also denote canonical inclusions into
the corresponding unital free products.

Now consider the homomorphisms 
\[ \bb\{ \pi_{t}^{T} \colon t \in S\} \colon
\bb\{ X_{t} \colon t \in S\} \to \bb\{ X_{t} \colon t \in T\} \]

\noindent and

\[ \bb\{ \pi_{t}^{T} \colon t \in T-S\} \colon
\bb\{ X_{t} \colon t \in T-S\} \to \bb\{ X_{t} \colon t \in T\} .\]

These two homomorphisms define the unique unital $\ast$-homomorphism

\[ f \colon \left(\bigstar_{{\mathbb C}}\{ X_{t}
\colon t \in S\}\right)
\bigstar_{{\mathbb C}}\left(\bigstar_{{\mathbb C}}\{
X_{t} \colon t \in T-S\}\right) \to \bb\{ X_{t} \colon t \in T\}\]

\noindent such that

\begin{equation}\label{EQ:1}
f\circ \pi_{1} = \bb\{ \pi_{t}^{T} \colon t \in S\}
\end{equation}

\noindent and

\begin{equation}\label{EQ:2}
f\circ \pi_{2} = \bb\{ \pi_{t}^{T} \colon t \in T-S\} .
\end{equation}
\medskip

\noindent Here $f = \left(\bb\{ \pi_{t}^{T} \colon t \in S\}\right)
\bb \left(\bb\{ \pi_{t}^{T} \colon t \in T-S\}\right)$.

Similarly consider the unique unital $\ast$-homomorphism
\[ g \colon \bb\{ X_{t} \colon t \in T\} \to
\left(\bigstar_{{\mathbb C}}\{ X_{t}
\colon t \in S\}\right)
\bigstar_{{\mathbb C}}\left(\bigstar_{{\mathbb C}}\{
X_{t} \colon t \in T-S\}\right)\]

\noindent satisfying the equalities

\begin{equation}\label{EQ:3}
g\circ \pi_{t}^{T} =
\begin{cases}
 \pi_{1}\circ \pi_{t}^{S}\;\;\; ,\;\text{if}\;\; t \in S ,\\
\pi_{2}\circ \pi_{t}^{T-S} ,\;\text{if}\;\; t \in T-S .
\end{cases}
\end{equation}

Next observe that if $t \in S$, then

\begin{equation}\label{EQ:4}
f \circ g \circ \pi_{t}^{T} \stackrel{(\ref{EQ:3})}{=}
f\circ\pi_{1}\circ \pi_{t}^{S} \stackrel{(\ref{EQ:1})}{=}
\bb\{ \pi_{t}^{T} \colon t \in S\} \circ \pi_{t}^{S} =
\pi_{t}^{T} .
\end{equation}

Similarly, if $t \in T-S$, then

\begin{equation}\label{EQ:5}
f \circ g \circ \pi_{t}^{T} \stackrel{(\ref{EQ:3})}{=}
f\circ\pi_{2}\circ \pi_{t}^{T-S} \stackrel{(\ref{EQ:2})}{=}
\bb\{ \pi_{t}^{T} \colon t \in T-S\} \circ \pi_{t}^{T-S} =
\pi_{t}^{T} .
\end{equation}

Now observe that (\ref{EQ:4}) and (\ref{EQ:5}) guarantee
the validity of the equality 

\begin{equation}\label{EQ:6}
f \circ g = \operatorname{id}_{\bb\{ X_{t} \colon t \in T\}} .
\end{equation}

In order to prove the equality 
\begin{equation}\label{EQ:7}
g\circ f = \operatorname{id}_{\left(\bb\{ \pi_{t}^{T} \colon t \in S\}\right)
\bb \left(\bb\{ \pi_{t}^{T} \colon t \in T-S\}\right)}
\end{equation}

\noindent it suffices to show that

\begin{equation}\label{EQ:8}
g \circ f \circ \pi_{1} = \pi_{1}
\end{equation}

\noindent and

\begin{equation}\label{EQ:9}
g \circ f \circ \pi_{2} = \pi_{2} .
\end{equation}

Note that (\ref{EQ:8}) follows from the following observation ($t \in S$): 

\begin{equation}\label{EQ:10}
g \circ f \circ \pi_{1}\circ \pi_{t}^{S}
\stackrel{(\ref{EQ:1})}{=}
g\circ \bb\{ \pi_{t}^{T} \colon t \in S\}\circ \pi_{t}^{S} =
g \circ \pi_{t}^{T} \stackrel{(\ref{EQ:3})}{=}
\pi_{1}\circ \pi_{t}^{S} .
\end{equation}

Similarly (\ref{EQ:9}) follows from the following
observation ($t \in T-S$):

\begin{equation}\label{EQ:11}
g \circ f \circ \pi_{2} \circ \pi_{t}^{T-S}
\stackrel{(\ref{EQ:2})}{=}
g\circ \bb\{ \pi_{t}^{T} \colon t \in S\}\circ \pi_{t}^{T-S} =
g \circ \pi_{t}^{T-S} \stackrel{(\ref{EQ:3})}{=}
\pi_{2}\circ \pi_{t}^{T-S} .
\end{equation}

\noindent This finishes proof of (\ref{EQ:7}). 

It only remains to note that, by (\ref{EQ:6}) and (\ref{EQ:7}),
both $f$ and $g$ are isomorphisms as required.
\end{proof}

\begin{lem}\label{L:1}
If $S \subseteq T$, then the unital
$\ast$-homomorphism

\[ \pi_{S}^{T} = \bb\{ \operatorname{id}_{X_{t}} \colon t \in S\}
\colon \bigstar_{{\mathbb C}}\{ X_{t}\colon t \in S\}
\to \bigstar_{{\mathbb C}}\{ X_{t}\colon t \in T\} \]
is injective.
\end{lem}
\begin{proof}
It can be shown, by applying the argument similar to the one used
in the proof of Lemma \ref{L:2}, that the homomorphism
$\pi_{S}^{T}$ coincides with the homomorphism

\[\pi_{1} \colon \bigstar_{{\mathbb C}}\{ X_{t}
\colon t \in S\} \to \left(\bigstar_{{\mathbb C}}\{ X_{t}
\colon t \in S\}\right)
\bigstar_{{\mathbb C}}\left(\bigstar_{{\mathbb C}}\{
X_{t} \colon t \in T-S\}\right) . \]
It only remains to note that $\pi_{1}$ is an inclusion by
\cite[Theorem 3.1]{blackadar1}.
\end{proof}

\begin{lem}\label{L:4}
If $\{ T_{\alpha} \colon \alpha < \tau\}$ is an increasing
well ordered collection of subsets of $T$ and
$T = \cup\{ T_{\alpha} \colon \alpha < \tau\}$, then
$\bigstar_{\mathbb C}\{ X_{t} \colon t \in T\}$ is
canonically isomorphic to
the direct limit of the well ordered direct system
$\left\{ \bigstar_{\mathbb C}\{ X_{t} \colon t \in T_{\alpha}\} ,
\pi_{T_{\alpha}}^{T_{\alpha +1}}, \tau \right\}$
\end{lem}
\begin{proof}
For each $\alpha <\tau$ consider the unital $\ast$-homomorphism
\[ \pi_{T_{\alpha}}^{T} \colon \bigstar_{\mathbb C}
\{ X_{t} \colon t \in T_{\alpha}\} \to
\bigstar_{\mathbb C}\{ X_{t} \colon t \in T\} ,\]
\noindent defined in Lemma \ref{L:1}.
Clearly
$\pi_{T_{\alpha +1}}^{T}\circ \pi_{T_{\alpha}}^{T_{\alpha +1}} =
\pi_{T_{\alpha}}^{T}$ for each $\alpha < \tau$. Consider
the unique unital $\ast$-homomorphism (see Subsection \ref{SS:morphisms})
\[ f \colon \varinjlim\left\{ \bigstar_{\mathbb C}
\{ X_{t} \colon t \in T_{\alpha}\} ,
\pi_{T_{\alpha}}^{T_{\alpha +1}}, \tau \right\} \to
\bigstar_{\mathbb C}\{ X_{t} \colon t \in T\}\]
\noindent such that $f \circ \pi_{\alpha} = \pi_{T_{\alpha}}$ for
each $\alpha < \tau$ (here 
\[ \pi_{\alpha} \colon \bigstar_{\mathbb C}\{ X_{t} \colon
t \in T_{\alpha}\} \to
\varinjlim\left\{ \bigstar_{\mathbb C}
\{ X_{t} \colon t \in T_{\alpha}\} ,
\pi_{ST_{\alpha}}^{T_{\alpha +1}}, \tau \right\}\]
denotes the $\alpha$-th limit injection of the above direct system).
Applying property $\bigl( \bb\bigr)$ it
is easy
to see that $f$ is an isomorphism.
\end{proof}

Finally we record the following statement. 

\begin{pro}\label{P:5}
Let
$\{X_{t} \colon t \in T\}$ be an infinite collection of  
unital $C^{\ast}$-algebras. Then the collection
$\left\{ \bigstar_{{\mathbb C}}\{ X_{t}
\colon t \in S\} , \pi_{S}^{R}, S, R
\in \exp_{<\omega}T \right\}$, consisting of the unital free products
of finite subcollections and above defined canonical injections,
is a direct system whose direct limit is the unital free product
$\bigstar_{{\mathbb C}}\{ X_{t} \colon t \in T\}$.

If the given collection $\{X_{t} \colon t \in T\}$ is uncountable
and consists of separable unital $C^{\ast}$-algebras, then the collection
$\left\{ \bigstar_{{\mathbb C}}\{ X_{t}
\colon t \in S\} , \pi_{S}^{R}, S, R
\in \exp_{\omega}T \right\}$, consisting of the unital free products
of countable subcollections and above defined canonical injections,
is a direct $C_{\omega}^{\ast}$-system
of the unital free product
$\bigstar_{{\mathbb C}}\{ X_{t} \colon t \in T\}$.
\end{pro}
\begin{proof}
The first part of this statement follows from the
above given definition of unital free products. In order to prove
the second part we need to show that
$\left\{ \bigstar_{{\mathbb C}}\{ X_{t}
\colon t \in S\} , \pi_{S}^{R}, S, R
\in \exp_{\omega}T \right\}$ is a direct $C_{\omega}^{\ast}$-system
associated with the unital free product
$\bigstar_{{\mathbb C}}\{ X_{t} \colon t \in T\}$. Let us verify
condition (a)--(d) of Definition \ref{D:smooth}. Condition (a)
is obvious since the set $\exp_{\omega}T$ is $\omega$-complete.
Conditions (c) and (d) follow from Lemmas \ref{L:1} and \ref{L:4}.
Finally condition (b), i.e. the fact that
$\bigstar_{{\mathbb C}}\{ X_{t} \colon t \in S\}$ is
separable for a countable subset $S \subseteq T$, follows
from \cite[Theorem 3.1]{blackadar1}.
\end{proof}

\begin{rem}\label{R:pointed}
The fact that the homomorphism $\pi_{S}^{T}$, indicated in Lemma \ref{L:1},
is injective can be significantly strengthened in the situation when each
$X_{t}$ admits a  unital $\ast$-homomorphism
$\varphi_{t} \colon X_{t} \to {\mathbb C}$. Indeed, in such a case, we can 
choose an index $s_{0} \in S$ and view the homomorphism $\varphi_{t}$ as a
unital $\ast$-homomorphism of $X_{t}$ into $X_{s_{0}}$. Next consider
the unital $\ast$-homomorphism
\[ g_{S}^{T} = \bb\{ g_{t} \colon t \in T\} \colon \bb\{ X_{t} \colon t \in T\}
\to \bb\{ X_{t} \colon t \in S\} ,\]

\noindent where
\[ g_{t} = 
\begin{cases}
\operatorname{id}_{X_{t}} \colon X_{t} \to X_{t} \; ,\; \text{if}\; t \in S ,\\
\varphi_{t} \colon X_{t} \to X_{s_{0}} \; ,\;\; \text{if}\; t \in T-S .
\end{cases}
\]
It is easy to show that
$g_{S}^{T}\circ \pi_{S}^{T} =
\operatorname{id}_{\bb\{ X_{t} \colon t \in S\}}$.
This means that $\pi_{S}^{T}$ is a coretraction and,
in particular, is injective.
\end{rem}

\begin{lem}\label{L:rem}
Let $X$ be a $C^{\ast}$-algebra admitting a unital
$\ast$-homomorphism into ${\mathbb C}$. If $Y$ is a $C^{\ast}$-subalgebra of $X$,
then $Y$ also admits a unital $\ast$-homomorphism into ${\mathbb C}$.
Projective unital
$C^{\ast}$-algebra admits a unital $\ast$-homomorphism into ${\mathbb C}$.
\end{lem}
\begin{proof}
The first part is trivial. If $X$ is a projective unital $C^{\ast}$-algebra,
then the projection $\pi_{1} \colon X \times {\mathbb C} \to X$ of
the direct product
$X\times {\mathbb C}$ onto the first coordinate has the
inverse, i.e. there exists
a unital $\ast$-homomorphism $i \colon X \to X \times {\mathbb C}$
such that
$\pi_{1}\circ i = \operatorname{id}_{X}$. Clearly the
projection $\pi_{2} \colon X\times {\mathbb C} \to {\mathbb C}$
onto the second coordinate is a unital $\ast$-homomorphism. It
only remains to note that the composition
$\pi_{2}\circ i \colon X \to {\mathbb C}$ is a unital
$\ast$-homomorphism.
\end{proof}


\bigskip
\bigskip

\section{Basic properties of doubly projective
homomorphisms and characterization of
projective unital $C^{\ast}$-algebras}\label{S:doubly}

Recall that a unital $C^{\ast}$-algebra $P$ is {\em projective} if
for any surjective unital $\ast$-homomorphism $p \colon X \to Y$
of unital $C^{\ast}$-algebras and
for any unital $\ast$-ho\-mo\-mor\-phism $f \colon P \to Y$ there exists
a unital $\ast$-homomorphism
$g \colon P \to X$ such that $p\circ g = f$. 


\subsection{Doubly projective homomorphisms}\label{SS:doubly}
The concept of doubly projective homomorphism was
introduced in \cite[Definition 3.1]{lope1}. In the definition
given below we do not assume that $X$ and $Y$ are
projective $C^{\ast}$-algebras.

\begin{defin}\label{D:doubly}
A unital $\ast$-homomorphism $i \colon X \to Y$ of unital
$C^{\ast}$-algebras $X$ and $Y$
is {\em doubly projective} if for any surjective unital 
$\ast$-ho\-mo\-mor\-phism
$p \colon A \to B$ between unital $C^{\ast}$-algebras
$A$ and $B$ and any two unital $\ast$-ho\-mo\-mor\-phisms
$f \colon X\to A$ and $g \colon Y \to B$ with $g\circ i = p\circ f$,
there exists a unital $\ast$-ho\-mo\-mor\-phism $h \colon B \to X$ such that
$f = h\circ i$ and $g = p\circ h$. In other words, any
commutative square diagram 

\[
\dgARROWLENGTH=2.0\dgARROWLENGTH
\begin{diagram}
\node{B}  \node{Y} \arrow{w,t}{g} \arrow{sw,t}{h}\\
\node{A} \arrow{n,l}{p} \node{X} \arrow{w,t}{f} \arrow{n,r}{i}
\end{diagram}
\]

\noindent with surjective $p$ can be completed by the diagonal arrow
with commuting triangles.
\end{defin}

We need some properties of doubly projective homomorphisms.

\begin{lem}\label{L:inj}
A doubly projective homomorphism $i \colon X \to Y$
of unital $C^{\ast}$-algebras
is a coretraction, i.e. there exists a unital $\ast$-homomorphism
$r \colon Y \to X$ such that
$r\circ i = \operatorname{id}_{X}$. In particular, a doubly projective homomorphism is injective.
\end{lem}
\begin{proof}
Let $i \colon X \to Y$ be a doubly projective homomorphism. Consider
the following commutative diagram

\[
\dgARROWLENGTH=2.0\dgARROWLENGTH
\begin{diagram}
\node{{\mathbf 0}}  \node{Y} \arrow{w,t}{g=
\operatorname{const}} \arrow{sw,t}{r}\\
\node{X} \arrow{n,l}{p=\operatorname{const}}
\node{X} \arrow{w,t}{f=\operatorname{id}_{X}} \arrow{n,r}{i}
\end{diagram}
\]

\noindent Since $i$ is doubly projective, there exists a
unital $\ast$-homomorphism
$r \colon Y \to X$ such that $r\circ i = \operatorname{id}_{X}$.
\end{proof}

\begin{lem}\label{L:iff}
Let $i \colon X \to Y$ be a doubly projective homomorphism of
unital $C^{\ast}$-algebras.
Then $X$ is projective if and only if $Y$ is projective.
\end{lem}
\begin{proof}
First suppose that $X$ is projective. In order to show that $Y$
is projective, consider a surjective unital $\ast$-homomorphism
$p \colon A \to B$
and a unital $\ast$-homomorphism $g \colon Y \to B$. Our goal is to find a
unital $\ast$-homomorphism
$\tilde{g} \colon Y \to A$ such that $p\circ \tilde{g} = g$. Since
$X$ is projective, there exists a unital $\ast$-homomorphism
$f \colon X \to A$ such that
$p\circ f = g\circ i$. Since $i$ is doubly projective there exists a
unital $\ast$-homomorphism $\tilde{g} \colon Y \to A$ such that
$g = p\circ \tilde{g}$
(and $f = \tilde{g}\circ i$). Obviously $\tilde{g}$ is a required
lift of $g$ and, consequently, $Y$ is projective.

Now assume that $Y$ is projective. In order to show that $X$ is
projective, consider a surjective unital $\ast$-homomorphism
$p \colon A \to B$ and a unital 
$\ast$-ho\-mo\-mor\-phism $f \colon X \to B$. Our goal is to find a
unital $\ast$-homomorphism
$\tilde{f} \colon X \to A$ such that $p\circ \tilde{f} = f$.
By Lemma \ref{L:inj}, there exists a unital $\ast$-homomorphism
$r \colon Y \to X$ such that $r\circ i = \operatorname{id}_{X}$.
Consider the composition $g = f \circ r \colon  Y \to B$.
Since $Y$ is projective, there exists a unital $\ast$-homomorphism
$\tilde{g} \colon Y \to A$
such that $p\circ \tilde{g} = g$. Let
$\tilde{f} = \tilde{g}\circ i$.
It only remains to note that $p\circ \tilde{f} =
p\circ \tilde{g} \circ i = g\circ i = f\circ r\circ i = f$.
\end{proof}

\begin{lem}\label{L:comp}
A finite composition of doubly projective homomorphisms is
doubly projective.
\end{lem}
\begin{proof}
Let $i_{1} \colon X_{1} \to X_{2}$ and
$i_{2} \colon X_{2} \to X_{3}$ be
doubly projective homomorphisms of unital $C^{\ast}$-algebras.
We need to show that
the composition $i = i_{2}\circ i_{1} \colon X_{1} \to X_{3}$
is also doubly
projective. Consider a surjective unital $\ast$-homomorphism
$p \colon A \to B$
and two unital $\ast$-homomorphisms
$g \colon X_{3} \to B$ and $f \colon X_{1} \to A$ such that
$g\circ i = p\circ f$. Consider the following commutative diagram

\[
\dgARROWLENGTH=2.0\dgARROWLENGTH
\begin{diagram}
\node{B}  \node{X_{2}} \arrow{w,t}{g\circ i_{2}} \arrow{sw,t}{f_{1}}\\
\node{A} \arrow{n,l}{p} \node{X_{1}} \arrow{w,t}{f} \arrow{n,r}{i_{1}}
\end{diagram}
\]

\noindent Since $i_{1}$ is doubly projective and since
$(g\circ i_{2})\circ i_{1} = g\circ i = p \circ f$, there
exists a unital $\ast$-homomorphism $f_{1} \colon X_{2} \to A$
such that
$p\circ f_{1} = g\circ i_{2}$ and $f = f_{1}\circ i_{1}$.

Next consider the commutative diagram

\[
\dgARROWLENGTH=2.0\dgARROWLENGTH
\begin{diagram}
\node{B}  \node{X_{3}} \arrow{w,t}{g} \arrow{sw,t}{h}\\
\node{A} \arrow{n,l}{p} \node{X_{2}} \arrow{w,t}{f_{1}}
\arrow{n,r}{i_{2}}
\end{diagram}
\]

\noindent Since $i_{2}$ is doubly projective and since
$p\circ f_{1} = g\circ i_{2}$, there exists a unital
$\ast$-homomorphism
$h \colon X_{3} \to A$ such that $p\circ h = g$ and
$f_{1} = h\circ i_{2}$.
It only remains to note that
\[ h\circ i = h\circ (i_{2}\circ i_{1}) = (h\circ i_{2})\circ i_{1}
= f_{1}\circ i_{1}= f .\]
\end{proof}

\begin{lem}\label{L:out}
Let $f \colon X \to Y$ be a doubly projective homomorphism.
Suppose that $f = f_{2}\circ f_{1}$, where $f_{2} \colon Z \to Y$
is a coretraction (i.e. there exists a unital $\ast$-homomorphism
$r \colon Y \to Z$ such that $r \circ f_{2} =
\operatorname{id}_{Z}$). Then $f_{1} \colon X \to Z$ is
also doubly projective.
\end{lem}
\begin{proof}
Let $p \colon A \to B$ be a surjective homomorphism
of unital $C^{\ast}$-algebras. Let also $g \colon X \to A$ and
$h \colon Z \to B$
be unital $\ast$-homomorphisms such that $p \circ g = h \circ f_{1}$. We need
to find a unital $\ast$-homomorphism $k \colon Z \to A$ such that
$k \circ f_{1} = g$ and $p \circ k = h$. Note that
\[ p\circ g = h \circ f_{1} = h \circ r \circ f_{2} \circ f_{1} =
(h \circ r) \circ f .\]
Since $f$ is doubly projective, there exists a unital $\ast$-homomorphism
$\widetilde{k} \colon Y \to A$ such that $\widetilde{k}\circ f = g$
and $p \circ \widetilde{k} = h \circ r$. Finally note that the composition
$k = \widetilde{k}\circ f_{2} \colon Z \to A$ has all the required properties.
Indeed,
\[ k\circ f_{1} = \widetilde{k}\circ f_{2}\circ f_{1} =
\widetilde{k}\circ f = g\]
\noindent and
\[ p\circ k = p\circ \widetilde{k}\circ f_{2} = h \circ r \circ f_{2} = h .\]
\end{proof}

\begin{lem}\label{L:retract}
Let $i \colon X \to Y$ be a unital $\ast$-homomorphism
which is a retract of a unital $\ast$-homomorphism
$i^{\prime} \colon X^{\prime} \to Y^{\prime}$. This means
that there exist
unital $\ast$-homomorphisms
$\varphi_{X} \colon X^{\prime} \to X$,
$\varphi_{Y} \colon Y^{\prime} \to Y$,
$\phi_{X} \colon X \to X^{\prime}$ and
$\phi_{Y} \colon Y \to Y^{\prime}$ such that
$i\circ \varphi_{X} = \varphi_{Y}\circ i^{\prime}$,
$i^{\prime}\circ\phi_{X} = \phi_{Y}\circ i$,
$\varphi_{X}\circ \phi_{X} = \operatorname{id}_{X}$ and
$\varphi_{Y}\circ \phi_{Y} = \operatorname{id}_{Y}$. In
other words the diagram

\[
\begin{CD}
Y @<\varphi_{Y}<< Y^{\prime} @<\phi_{Y}<< Y\\
@A{i}AA @AA{i^{\prime}}A @AA{i}A\\
X @<\varphi_{X}<< X^{\prime} @<\phi_{X}<< X
\end{CD}
\]

\noindent commutes. In this situation, if $i^{\prime}$ is
doubly projective, then
$i$ is also doubly projective.
\end{lem}
\begin{proof}
Consider a surjective unital $\ast$-homomorphism $p \colon A \to B$
and two unital $\ast$-homomorphisms $f \colon X \to A$ and $g \colon Y \to B$
such that $g \circ i = p\circ f$. Here is the corresponding diagram

\[
\begin{CD}
B @<g<< Y @<\varphi_{Y}<< Y^{\prime} @<\phi_{Y}<< Y\\
@A{p}AA @A{i}AA @AA{i^{\prime}}A @AA{i}A\\
A @<f<< X @<\varphi_{X}<< X^{\prime} @<\phi_{X}<< X
\end{CD}
\]

\noindent Let $f^{\prime} = f \circ \varphi_{X} \colon X^{\prime} \to A$
and $g^{\prime} = g\circ \varphi_{Y} \colon Y^{\prime} \to B$. Note that

\[ g^{\prime}\circ i^{\prime} = g \circ \varphi_{Y}\circ i^{\prime}
= g\circ i \circ \varphi_{X} =
p\circ f \circ \varphi_{X} = p\circ f^{\prime}.\]
Since $i^{\prime}$ is doubly projective, there exists a
unital $\ast$-homomorphism
$h^{\prime} \colon Y^{\prime} \to A$ such that
$p\circ h^{\prime} = g^{\prime}$ and
$h^{\prime}\circ i^{\prime} = f^{\prime}$. 

Now consider the composition
$h = h^{\prime}\circ\phi_{Y} \colon Y \to B$ and observe that

\[ p \circ h = p\circ h^{\prime}\circ \phi_{Y} = g^{\prime}\circ \phi_{Y} = g\circ\varphi_{Y}\circ \phi_{Y} = g \]

\noindent and

\[ h \circ i = h^{\prime}\circ \phi_{Y} \circ i =
h^{\prime}\circ i^{\prime}\circ \phi_{X} = f^{\prime}\circ \phi_{X} = f\circ\varphi_{X}\circ\phi_{X} = f .\]
This shows that $i$ is doubly projective.
\end{proof}

The following statement provides an important class of doubly projective homomorphisms.

\begin{lem}\label{L:ex}
Let $X$ be a unital $C^{\ast}$-algebra and $Y$ be a projective unital
$C^{\ast}$-algebra.
Then the canonical inclusion
$\pi_{X} \colon X \hookrightarrow X \bigstar_{{\mathbb C}} Y$
is doubly projective.
\end{lem}
\begin{proof}
Consider a surjective unital $\ast$-ho\-mo\-mor\-phism $p \colon A \to B$ and
two unital $\ast$-ho\-mo\-mor\-phisms
$f \colon X \to A$ and $g \colon X \bigstar_{{\mathbb C}} Y \to B$ such that
$p\circ f = g \circ i_{X}$. Our goal is to construct a unital $\ast$-homomorphism
$h \colon X \bigstar_{{\mathbb C}} Y \to A$ such that $p\circ h = g$ and
$h\circ \pi_{X} = f$. Let $\pi_{Y} \colon Y \to X \bigstar_{{\mathbb C}} Y$
denote the
canonical embedding of $Y$ into $X \bigstar_{{\mathbb C}} Y$.
Since $Y$ is projective,
there exists a unital $\ast$-homomorphism $h_{1} \colon Y \to A$ such that
$p\circ h_{1} = g\circ i_{Y}$. The two unital $\ast$-homomorphisms
$f \colon X \to A$
and $h_{1} \colon Y \to A$ define the unique unital $\ast$-homomorphism
$h \colon X \bigstar_{{\mathbb C}} Y \to A$
such that $h\circ \pi_{X} = f$ and $h\circ \pi_{Y} = h_{1}$.
Finally, observe that
$g\circ \pi_{X} = p\circ f = (p\circ h)\circ \pi_{X}$ and
$g\circ \pi_{Y} = p\circ h_{1} = (p \circ h)\circ \pi_{Y}$. This shows that
$p\circ h = g$.
\end{proof}

Next we introduce the concept of a doubly projective
homomorphism of separable type.

\begin{defin}\label{D:septype}
We say that a doubly projective homomorphism
$i \colon X \to Y$ between projective unital $C^{\ast}$-algebras
has a separable type, if there exist a
projective unital $C^{\ast}$-algebra $X^{\prime}$ such that
$d(X^{\prime}) = d(X)$, a separable
projective unital $C^{\ast}$-algebra $Y^{\prime}$ and two
surjective unital $\ast$-homomorphisms
$\varphi_{X} \colon X^{\prime} \to X$ and
$\varphi_{Y} \colon X^{\prime}\bb Y^{\prime} \to Y$ such that
$i \circ \varphi_{X} = \varphi_{Y}\circ \pi_{X^{\prime}}$, where
$\pi_{X^{\prime}} \colon X^{\prime} \to  X^{\prime} \bb Y^{\prime}$
denotes the natural inclusion. In other words we require
the commutativity of the following diagram

\[
\begin{CD}
Y @<\varphi_{Y}<<
X^{\prime}\bb Y^{\prime}\\
@A{i}AA @AA{\pi_{X^{\prime}}}A \\
X @<\varphi_{X}<<
X^{\prime}.
\end{CD}
\]
\end{defin}

\begin{lem}\label{L:sepsep}
Every doubly projective homomorphism between
separable projective unital
$C^{\ast}$-algebras has a separable type.
\end{lem}
\begin{proof}
Let $i \colon X \to Y$ be a doubly projective homomorphism and
$X$ and $Y$ be separable projective unital $C^{\ast}$-algebras.
Consider the unital free product $X \bb Y$ and note that the diagram

\[
\begin{CD}
Y @<i\bb\operatorname{id}_{Y}<<
X\bb Y\\
@A{i}AA @AA{\pi_{X}}A \\
X @<\operatorname{id}_{X}<<
X
\end{CD}
\]

\noindent commutes. Also observe that $X\bb Y$ is a
projective unital $C^{\ast}$-algebra. Clearly
$i\bb \operatorname{id}_{Y} \colon X\bb Y \to X$ is
surjective, because 
$(i\bb \operatorname{id}_{Y})\circ \pi_{Y} = \operatorname{id}_{Y}$
is surjective. 
\end{proof}

\begin{lem}\label{L:well}
Let ${\mathcal S} = \{ X_{\alpha}, i_{\alpha}^{\alpha +1},\tau\}$ be a
well-ordered continuous direct system of unital $C^{\ast}$-algebras. If the
short injection
$i_{\alpha}^{\alpha +1} \colon X_{\alpha} \to X_{\alpha +1}$ of the
system $\mathcal S$ is doubly projective
for each $\alpha < \tau$, then the limit injection
$i_{0} \colon X_{0} \to \varinjlim{\mathcal S}$ is also doubly projective.
\end{lem}
\begin{proof}
Let $p \colon A \to B$ be a surjective unital $\ast$-homomorphism of
unital $C^{\ast}$-algebras. Let also $g \colon X_{0} \to A$ and
$h \colon \varinjlim{\mathcal S}_{X} \to B$ be unital $\ast$-homomorphisms
such that $p \circ g = h \circ i_{0}$. 

By induction we construct a well ordered collection
$\{ k_{\alpha} \colon X_{\alpha} \to A ; \alpha < \tau \}$ of
unital $\ast$-homomorphisms.
Let $k_{0} = g$ and suppose that we have already
constructed $\ast$-homomorphisms
$k_{\alpha}$ for each $\alpha < \gamma$, where $\gamma < \tau$,
in such a way that the following conditions are satisfied:
\begin{itemize}
\item[(a)]
$k_{\alpha} = k_{\alpha +1}\circ i_{\alpha}^{\alpha +1}$
for each $\alpha < \gamma$.
\item[(b)]
$p\circ k_{\alpha} = h \circ i_{\alpha}$ for each $\alpha < \tau$.
\item[(c)]
$k_{\beta} = \varinjlim\{ k_{\alpha}; \alpha < \beta \}$
whenever $\beta$ is a limit
ordinal number with $\beta < \gamma$.
\end{itemize}

Let us construct a $\ast$-homomorphism $k_{\gamma} \colon X_{\gamma} \to A$.

If $\gamma$ is a limit ordinal number, then
let $k_{\gamma} = \varinjlim\{ k_{\alpha}; \alpha < \gamma\}$.

If $\gamma = \alpha +1$, then consider the following commutative diagram

\[
\begin{CD}
A @>p>> B\\
@A{k_{\alpha}}AA @AA{h\circ i_{\alpha}}A\\
X_{\alpha} @>i_{\alpha}^{\alpha +1}>> X_{\alpha +1}
\end{CD}
\]

\noindent Since $i_{\alpha}^{\alpha +1} \colon
X_{\alpha} \to X_{\alpha +1}$ is doubly projective there exists
a unital $\ast$-homomorphism $k_{\alpha +1} \colon X_{\alpha +1} \to A$
such that $k_{\alpha} = k_{\alpha +1}\circ i{\alpha}^{\alpha +1}$ and
$p\circ k_{\alpha +1} = h \circ i_{\alpha}$.

Thus, the homomorphisms $k_{\alpha} \colon X_{\alpha} \to A$ are
constructed for each $\alpha < \tau$ and satisfy the above
stated properties for each $\alpha < \tau$. It only remains to note that
for the unital $\ast$-homomorphism
$k = \varinjlim\{ k_{\alpha} ;\alpha < \tau\} \colon
\varinjlim{\mathcal S}_{X} \to A$ we have
$g = k_{0} = k\circ i_{0}$ and $h = p\circ k$
as required.
\end{proof}

As was pointed out in the introduction, there is a deeper relation between
doubly projective homomorphisms and projective $C^{\ast}$-algebras, than it
might appear to be the case. Let $\operatorname{Mor}({\mathcal C}_{1}^{\ast})$ denote
the category of unital $\ast$-homomorphisms of unital $C^{\ast}$-algebras.
The following statement is true.

\begin{pro}\label{P:morphisms}
The following conditions are equivalent for a unital
$\ast$-ho\-mo\-mor\-phism
$f \colon X \to Y$ of projective unital $C^{\ast}$-algebras:
\begin{itemize}
\item[(a)]
$f$ is doubly projective.
\item[(b)]
$f$ is a projective object of the category
$\operatorname{Mor}({\mathcal C}_{1}^{\ast})$.
\end{itemize}
\end{pro}
\begin{proof}
(a) $\Longrightarrow$ (b). Let $p \colon A\to B$ and
$q \colon C \to D$ be objects of the category
$\operatorname{Mor}({\mathcal C}_{1}^{\ast})$ and $(s,r) \colon p \to q$ be
an epimorphism of the same category. Our goal is to show that
for any morphism $(\alpha ,\beta ) \colon f \to q$ of
$\operatorname{Mor}({\mathcal C}_{1}^{\ast})$ there exists a morphism
$(\widetilde{\alpha},\widetilde{\beta}) \colon f \to p$ of
$\operatorname{Mor}({\mathcal C}_{1}^{\ast})$ such that
$(s,r)\circ (\widetilde{\alpha},\widetilde{\beta}) = (\alpha ,\beta )$.

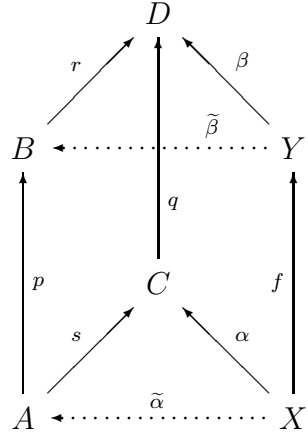
\[
\dgARROWLENGTH=2.0\dgARROWLENGTH
\begin{diagram}
\node[2]{D} \\
\node{B} \arrow{ne,t}{r} \node[2]{Y} \arrow{nw,t}{\beta}
\arrow[2]{w,t,..,1}{\widetilde{\beta}} \\
\node[2]{C} \arrow[2]{n,r,1}{q}\\
\node{A} \arrow[2]{n,r}{p} \arrow{ne,t}{s} \node[2]{X} \arrow[2]{n,l}{f} \arrow[2]{w,t,..}{\widetilde{\alpha}} \arrow{nw,t}{\alpha}
\end{diagram}
\]

\noindent Since $(r,s)$ is an epimorphism in
$\operatorname{Mor}({\mathcal C}_{1}^{\ast})$ it follows that each
of the homomorphisms $r$ and $s$ is surjective. Since $X$ is projective,
there exists a unital
$\ast$-homomorphism $\widetilde{\alpha} \colon X \to A$
such that $s \circ\widetilde{\alpha} = \alpha$. Clearly
$r\circ p \circ \widetilde{\alpha} =
q\circ s \circ \widetilde{\alpha} = q\circ \alpha = \beta \circ f$.
Consequently, since $f$ is doubly projective, there exists a
unital $\ast$-homomorphism $\widetilde{\beta} \colon Y \to B$ such that 
$\widetilde{\beta}\circ f = p \circ\widetilde{\alpha}$ and
$r \circ \widetilde{\beta} = \beta$. In other words the following diagram

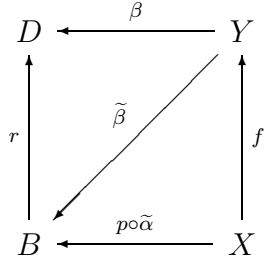
\[
\dgARROWLENGTH=2.0\dgARROWLENGTH
\begin{diagram}
\node{D}  \node{Y} \arrow{w,t}{\beta} \arrow{sw,t}{\widetilde{\beta}}\\
\node{B} \arrow{n,l}{r} \node{X}
\arrow{w,t}{p\circ\widetilde{\alpha}} \arrow{n,r}{f}
\end{diagram}
\]

\noindent commutes. The straitforward verification shows that
$(s,r)\circ (\widetilde{\alpha},\widetilde{\beta}) = (\alpha ,\beta )$
as required.

(b) $\Longrightarrow$ (a). Now suppose that $f \colon X \to Y$ is a projective
object of the category $\operatorname{Mor}({\mathcal C}_{1}^{\ast})$.
In order to show that $f$ is doubly projective, consider a surjective
unital $\ast$-homomorphism $p \colon A \to B$ and two unital
$\ast$-homomorphisms $g \colon X \to A$ and $h \colon Y \to B$ such that
$p\circ g = h\circ f$. Clearly the pair $(g,h)$ forms a
morphism $(g,h) \colon f \to p$ in the category
$\operatorname{Mor}({\mathcal C}_{1}^{\ast})$. Consider also
the epimorphism (in the category $\operatorname{Mor}({\mathcal C}_{1}^{\ast})$)
$(\operatorname{id}_{A}, p) \colon \operatorname{id}_{A} \to p$ sending the left
vertical arrow in the following diagram onto the middle one.

\[
\begin{CD}
A @>p>> B @<h<< Y\\
@A{\operatorname{id}_{A}}AA @AA{p}A @AA{f}A\\
A @>\operatorname{id}_{A}>> A @<g<<X
\end{CD}
\]

\noindent Since $f$ is a projective object in the category
$\operatorname{Mor}({\mathcal C}_{1}^{\ast})$, it follows that there exists a
morphism
$(\widetilde{g},\widetilde{h}) \colon f \to \operatorname{id}_{A}$,
consisting of the unital $\ast$-homomorphisms
$\widetilde{g} \colon X \to A$ and $\widetilde{h} \colon Y \to A$,
such that
$(\operatorname{id}_{A}, p)\circ (\widetilde{g},\widetilde{h}) = (g,h)$.
This implies that $\operatorname{id}_{A}\circ \widetilde{g} = g$, i.e.
$\widetilde{g} = g$, and $p \circ \widetilde{h} = h$.
In order to prove the equality
$\widetilde{h}\circ f = g$, simply note that
$(\widetilde{g},\widetilde{h}) \colon f \to
\operatorname{id}_{A}$ is a morphism
in the category $\operatorname{Mor}({\mathcal C}_{1}^{\ast})$.
Thus $f$ is doubly projective.
\end{proof}


\subsection{Characterization of projective unital
$C^{\ast}$-algebras}\label{SS:charactpa}
We begin with the following preliminary result.
\begin{lem}\label{L:first}
Let $X$ be a projective unital $C^{\ast}$-algebra of
density $\tau > \omega$. Then $X$ admits a direct
$C_{\omega}^{\ast}$-system
${\mathcal S}_{X} = \{ X_{t}, i_{s}^{t}, {\mathcal A}\}$,
consisting of separable projective unital 
$C^{\ast}$-subalgebras of $X$. We may assume that ${\mathcal A}$
is a cofinal and $\omega$-closed subset of $\exp_{\omega}\tau$.
\end{lem}
\begin{proof}
Let $A$ be a dense subset
of $X$ such that $|A| = \tau$. Let also $T = \exp_{\omega}A$.
Since $\tau > \omega$, it follows that $|T| = \tau$.
As in the proof of
Proposition \ref{P:exists}, we can conclude that $X$ is the limit of the
direct system $\{ X_{t}, i_{s}^{t}, T\}$,
consisting of separable unital
$C^{\ast}$-subalgebras of $X$ (generated by countable subsets of $A$)
and associated inclusion maps.

Next consider the unital $\ast$-homomorphism
$\varphi \colon \bigstar_{{\mathbb C}}\{ X_{t} \colon t \in T\} \to X$,
generated
by the homomorphisms $i_{t} \colon X_{t} \to X$. This means
that $\varphi\circ \pi_{X_{t}} = i_{t}$ for each $t \in T$
(here $\pi_{X_{t}} \colon X_{t} \to
\bigstar_{{\mathbb C}}\{ X_{t} \colon t \in T\}$
denotes the canonical inclusion). Note that $\varphi$ is a
surjective unital $\ast$-homomorphism. This follows
from Lemma \ref{L:strong}.

Recall that by Proposition \ref{P:5}, the collection 

\[ {\mathcal S} = \left\{ \bigstar_{{\mathbb C}}\{ X_{t}
\colon t \in S\} , \pi_{S}^{R}, \exp_{\omega}T \right\} \]

\noindent is a direct $C_{\omega}^{\ast}$-system such that
$\bigstar_{{\mathbb C}}\{ X_{t} \colon t \in T\} = \varinjlim{\mathcal S}$.

For each $S \subseteq T$ let $X_{S} =
\operatorname{cl}_{X}\varphi\left(\bigstar_{{\mathbb C}}\{ X_{t}
\colon t \in S\}\right)$. Also by
$\varphi_{S} \colon \bigstar_{{\mathbb C}}\{ X_{t}
\colon t \in S\} \to X_{S}$ we denote the restriction of
the homomorphism $\varphi$ onto the unital free product
$\bigstar_{{\mathbb C}}\{ X_{t}
\colon t \in S\}$. We have the following commutative diagram

\[
\begin{CD}
X @<\varphi<<
\bigstar_{{\mathbb C}}\{ X_{t} \colon t \in T\}\\
@A{i_{S}}AA @AA{\pi_{S}^{T}}A \\
X_{S} @<\varphi_{S}<<
\bigstar_{{\mathbb C}}\{ X_{t} \colon t \in S\} ,
\end{CD}
\]

\noindent where $i_{S} \colon X_{S} \to X$ denotes the inclusion.

It is obvious that the system
${\mathcal S}_{X} = \{ X_{S}, i_{S}^{R}, \exp_{\omega}T\}$, consisting of
$C^{\ast}$-sub\-algebras
$X_{S}$ of $X$ and their natural inclusions
$i_{S}^{R} \colon X_{S} \to X_{R}$, forms a
direct $C^{\ast}_{\omega}$-system such that
$X = \varinjlim{\mathcal S}_{X}$. Also note that
\[ \{ \varphi_{S} \colon S \in \exp_{\omega}T\}
\colon {\mathcal S} \to {\mathcal S}_{X}\]

\noindent is a morphism between the indicated direct
systems such that
$\varphi = \varinjlim\{ \varphi_{S} \colon S \in \exp_{\omega}T\}$.

Since $X$ is a projective
$C^{\ast}$-algebra, there exists a unital $\ast$-ho\-mo\-mor\-phism
$\phi \colon X\to \bigstar_{{\mathbb C}}\{ X_{t} \colon t \in T\}$
such that
$\varphi \circ \phi = \operatorname{id}_{X}$. 

According to Theorem \ref{T:spectral}, applied to the homomorphism
$\phi \colon \varinjlim{\mathcal S}_{X} \to \varinjlim{\mathcal S}$, 
there exist a cofinal and $\omega$-closed subset ${\mathcal A}$ of
$\exp_{\omega}\tau$ and a morphism
\[ \{ \phi_{S} \colon S \in {\mathcal A}\}
\colon {\mathcal S}_{X}|{\mathcal A} \to {\mathcal S}|{\mathcal A} \]

\noindent such that
$\phi = \varinjlim\{ \phi_{S} \colon S \in {\mathcal A}\}$.
In particular, the square diagram

\[
\begin{CD}
X @>\phi>>
\bigstar_{{\mathbb C}}\{ X_{t} \colon t \in T\}\\
@A{i_{S}}AA @AA{\pi_{S}^{T}}A \\
X_{S} @>\phi_{S}>>
\bigstar_{{\mathbb C}}\{ X_{t} \colon t \in S\} ,
\end{CD}
\]

\noindent commutes for each $S \in {\mathcal A}$. Note also that
$\varphi_{S}\circ \phi_{S} = \operatorname{id}_{X_{S}}$ for
each $S \in {\mathcal A}$.

According to Lemma \ref{L:rem}, the $C^{\ast}$-algebra $X$,
and hence each $X_{t}$,
$t \in T$, admits a unital $\ast$-homomorphism into ${\mathbb C}$.
Consequently, by Remark \ref{R:pointed}, the inclusion
$\pi_{S}^{T} \colon \bb\{ X_{t} \colon t \in S\} \to
\bb\{ X_{t} \colon t \in T\}$ is a coretraction with the
associated retraction $g_{S}^{T} \colon
\bb\{ X_{t} \colon t \in T\} \to \bb\{ X_{t} \colon t \in S\}$.
Consider the unital $\ast$-homomorphism $r_{S} \colon X \to X_{S}$,
defined as the composition $r_{S} = \varphi_{S}\circ g_{S}^{T}\circ \phi$.
Note that
\[ r_{S} \circ i_{S} = \varphi_{S}\circ g_{S}^{T}\circ \phi \circ i_{S} = 
\varphi_{S}\circ g_{S}^{T}\circ \pi_{S}^{T}\circ \phi_{S} =
\varphi_{S}\circ \phi_{S} = \operatorname{id}_{X_{S}} ,\]

\noindent which shows that $r_{S}$ is a retraction. It only
remains to note that $X_{S}$, as a retract of $X$, is projective. 
\end{proof}

The following statement provides a characterization of non-separable
projective unital $C^{\ast}$-algebras. It should be noted
that condition (b) of Theorem \ref{T:charact} is significantly
stronger than the conclusion of Lemma \ref{L:first}.

\begin{thm}\label{T:charact}
The following conditions are equivalent for any unital $C^{\ast}$-algebra
$X$ of density $d(X) = \tau > \omega$:
\begin{itemize}
\item[(a)]
$X$ is projective.
\item[(b)]
$X$ is isomorphic to the limit of a direct $C_{\omega}^{\ast}$-system
${\mathcal S}_{X} = \{ X_{\alpha}, i_{\alpha}^{\beta}, {\mathcal A}\}$,
consisting of
separable projective unital $C^{\ast}$-algebras $X_{\alpha}$
and doubly projective limit injections
$i_{\alpha} \colon X_{\alpha} \to X$, $\alpha \in {\mathcal A}$.
We may assume that ${\mathcal A}$ is cofinal and $\omega$-closed in
$\exp_{\omega}\tau$.
\item[(c)]
$X$ is isomorphic to the limit of a well-ordered continuous
direct system ${\mathcal S}_{X} = \{ X_{\alpha},
i_{\alpha}^{\alpha+1}, \tau\}$ of length
$\tau$ satisfying the following properties:
\begin{enumerate}
\item
$X_{\alpha}$ is a projective unital $C^{\ast}$-algebra for each $\alpha < \tau$.
\item
Short injection $i_{\alpha}^{\alpha +1} \colon X_{\alpha} \to X_{\alpha +1}$
is doubly projective and has a separable type for each $\alpha < \tau$.
\item
$X_{0}$ is a separable projective unital $C^{\ast}$-algebra.
\end{enumerate}
\end{itemize}
\end{thm}
\begin{proof}
{\em Part I}. First we show that if $X$ is a projective
$C^{\ast}$-algebra,
then there exists a well ordered continuous direct system
${\mathcal S}_{X} = \{ X_{\alpha}, i_{\alpha}^{\alpha +1}, \tau\}$,
satisfying condition (c). While proving this we will show the existence
of a direct $C_{\omega}^{\ast}$-system
$\{ X_{\alpha}, i_{\alpha}^{\beta}, {\mathcal A}\}$,
satisfying condition (b).

According to condition (a) and Lemmas \ref{L:first} and
\ref{L:strong} there exists a collection
$\{ X_{t} \colon t \in T\}$, consisting of separable
unital projective
$C^{\ast}$-subalgebra of $X$, such that
$X = \cup\{ X_{t} \colon t \in T\}$ and $|T| = \tau$. 

Below we follow the proof of Lemma \ref{L:first}. The fact that each 
$X_{t}$, $t \in {\mathcal A}$, is projective becomes crucial
later in this proof.

As in the proof of Lemma \ref{L:first}, the homomorphisms
$i_{t} \colon X_{t} \to X$, $t \in T$, generate the surjective unital
$\ast$-homomorphism
$\varphi \colon \bigstar_{{\mathbb C}}\{ X_{t} \colon t \in T\} \to X$
such that $\varphi\circ \pi_{X_{t}} = i_{t}$ for each $t \in T$
(here $\pi_{X_{t}} \colon X_{t} \to
\bigstar_{{\mathbb C}}\{ X_{t} \colon t \in T\}$
denotes the canonical inclusion).

Recall that by Proposition \ref{P:5}, the collection 

\[ {\mathcal S} = \left\{ \bigstar_{{\mathbb C}}\{ X_{t}
\colon t \in S\} , \pi_{S}^{R}, \exp_{\omega}T \right\} \]

\noindent is a direct $C_{\omega}^{\ast}$-system such that
$\bigstar_{{\mathbb C}}\{ X_{t} \colon t \in T\} = \varinjlim{\mathcal S}$.

For each $S \subseteq T$ let $X_{S} =
\operatorname{cl}_{X}\varphi\left(\bigstar_{{\mathbb C}}\{ X_{t}
\colon t \in S\}\right)$. Also by
$\varphi_{S} \colon \bigstar_{{\mathbb C}}\{ X_{t}
\colon t \in S\} \to X_{S}$ we denote the restriction of
the homomorphism $\varphi$ onto the unital free product
$\bigstar_{{\mathbb C}}\{ X_{t}
\colon t \in S\}$. We have the following commutative diagram

\[
\begin{CD}
X @<\varphi<<
\bigstar_{{\mathbb C}}\{ X_{t} \colon t \in T\}\\
@A{i_{S}}AA @AA{\pi_{S}^{T}}A \\
X_{S} @<\varphi_{S}<<
\bigstar_{{\mathbb C}}\{ X_{t} \colon t \in S\} ,
\end{CD}
\]

\noindent where $i_{S} \colon X_{S} \to X$ denotes the inclusion.

It is obvious that the system
${\mathcal S}_{X} = \{ X_{S}, i_{S}^{R},\exp_{\omega}T\}$, consisting of
$C^{\ast}$-sub\-algebras
$X_{S}$ of $X$ and their natural inclusions
$i_{S}^{R} \colon X_{S} \to X_{R}$, forms a
direct $C^{\ast}_{\omega}$-system such that $X = \varinjlim{\mathcal S}_{X}$.

Since, by (a), $X$ is a projective
$C^{\ast}$-algebra, there exists a unital $\ast$-ho\-mo\-mor\-phism
$\phi \colon X\to \bigstar_{{\mathbb C}}\{ X_{t} \colon t \in T\}$
such that
$\varphi \circ \phi = \operatorname{id}_{X}$. 

Let us say that a subset $S \subseteq T$ is admissible if
$\phi (X_{S}) \subseteq \bb\{ X_{t} \colon t \in S\}$.
This clearly means that the diagram

\[
\begin{CD}
X @>\phi>>
\bigstar_{{\mathbb C}}\{ X_{t} \colon t \in T\}\\
@A{i_{S}}AA @AA{\pi_{S}^{T}}A \\
X_{S} @>\phi_{S}>>
\bigstar_{{\mathbb C}}\{ X_{t} \colon t \in S\} ,
\end{CD}
\]

\noindent where
$\phi_{S} = \phi |X_{S} \colon X_{S} \to
\bb\{ X_{t} \colon t \in S\}$, commutes.
\medskip

We need to state some of the properties of admissible subsets.

\medskip

{\bf Claim 1.} {\em If $S \subseteq T$ is admissible,
then $\phi_{S}\circ \varphi_{S} = \operatorname{id}_{X_{S}}$.}

{\em Proof of Claim 1.} Follows form the above constructions
and the equality
$\phi\circ\varphi = \operatorname{id}_{X}$ (see the proof
of Lemma \ref{L:first}).

{\bf Claim 2.} {\em If $S$ is an admissible subset of
 $T$, then $X_{S} = \varphi\left(\bb\{ X_{t} \colon t \in S\}\right)$.}

{\em Proof of Claim 2.} Follows from Claim 1 (see the proof
of Lemma \ref{L:first}).

{\bf Claim 3.} {\em The union of an arbitrary collection of
admissible subsets of $T$ is admissible.} 

{\em Proof of Claim 3.} Let $S_{i}$, $i \in I$ be an admissible
subset of $T$ and let $S = \bigcup\{ S_{i} \colon i \in I\}$.
First observe that
\begin{multline*}
X_{S} = \operatorname{cl}_{X}\varphi\left(\bb\{ X_{t} \colon t \in S\}\right)
= \operatorname{cl}_{X}\varphi\left( \bigcup\{
\bb\{ X_{t} \colon t \in S_{i}\} \colon i \in I\}\right) =\\
 \operatorname{cl}_{X}\left( \bigcup\{ \varphi\left(
\{ X_{t} \colon t \in S_{i}\}\right) \colon i \in I\right)
\subseteq \operatorname{cl}_{X}\left( \bigcup\{
X_{S_{i}} \colon i \in I\}\right) .
\end{multline*}

\noindent Consequently

\begin{multline*}
\phi (X_{S}) \subseteq \phi\left(\operatorname{cl}_{X}
\left( \bigcup\{ X_{S_{i}} \colon i \in I\}\right)\right)
\subseteq \operatorname{cl}_{\bb\{ X_{t} \colon t \in T\}}
\left( \phi\left(\bigcup\{ X_{S_{i}} \colon i \in I\}\right)\right)
\subseteq\\
\operatorname{cl}_{\bb\{ X_{t} \colon t \in T\}}\left(
\bigcup\{ \phi (X_{S_{i}} \colon i \in I\}\right)
\subseteq \operatorname{cl}_{\bb\{ X_{t} \colon t \in T\}}
\left(\bigcup\{ \bb\{ X_{t} \colon t \in S_{i}\} \colon
i \in I\}\right)\\ \subseteq
 \bb\{ X_{t} \colon t \in S\} .
\end{multline*}

{\bf Claim 4.} {\em If $S$ is an admissible subset of $T$,
then $X_{S}$ is a projective $C^{\ast}$-subalgebra of $X$.}

{\em Proof of Claim 4.} See the proof of lemma \ref{L:first}.

{\bf Claim 5.} {\em Every countable subset of $T$ is contained
in a countable admissible subset of $T$.}

{\em Proof of Claim 5.} According to Theorem \ref{T:spectral},
applied to the homomorphism
$\phi \colon \varinjlim{\mathcal S}_{X} \to \varinjlim{\mathcal S}$, 
there exist a cofinal and $\omega$-closed subset ${\mathcal A}$ of
$\exp_{\omega}\tau$ and a morphism
\[ \{ \phi_{S} \colon S \in {\mathcal A}\}
\colon {\mathcal S}_{X}|{\mathcal A} \to {\mathcal S}|{\mathcal A} \]

\noindent such that
$\phi = \varinjlim\{ \phi_{S} \colon S \in {\mathcal A}\}$.
Clearly each $S \in {\mathcal A}$ is admissible.

{\bf Claim 6.} {\em If $S$ is an admissible subset of
$T$, then the inclusion
$i_{S} \colon X_{S} \to X$ is doubly projective.}

{\em Proof of Claim 6.} Recall that the following diagram

\[
\begin{CD}
X @<\varphi<<
\bigstar_{{\mathbb C}}\{ X_{t} \colon t \in T\} @<\phi<< X\\
@A{i_{S}}AA @AA{\pi_{S}^{T}}A @AA{i_{S}}A \\
X_{S} @<\varphi_{S}<<
\bigstar_{{\mathbb C}}\{ X_{t} \colon t \in S\} @<\phi_{S}<< X_{S}
\end{CD}
\]

\noindent commutes and that $\varphi\circ \phi = \operatorname{id}_{X}$
and $\varphi_{S}\circ\phi_{S} = \operatorname{id}_{X_{S}}$.

Since each $X_{t}$, $t \in T$, is projective (this is where
Lemma \ref{L:first} is actually being used) we easily
conclude that $\bb \{ X_{t} \colon t \in T-S\}$ is also
projective (compare to \cite[Propositions 2.31, 2.32]{blackadar2}).
By Lemmas \ref{L:2} and \ref{L:ex}, the inclusion
$\pi_{S}^{T} \colon \bb\{ X_{t} \colon t \in S\} \to
\bb\{ X_{t} \colon t \in T\}$ is doubly projective. Finally,
Lemma \ref{L:retract} guarantees that the inclusion
$i_{S} \colon X_{S} \to X$ is also doubly projective.
This completes proof of Claim 6.

Now consider the direct system
${\mathcal S}_{X}|{\mathcal A} = \{ X_{S}, i_{S}^{R}, {\mathcal A}\}$.
Clearly ${\mathcal S}_{X}|{\mathcal A}$ is a
direct $C_{\omega}^{\ast}$-system such that
$X = \varinjlim{\mathcal S}_{X}|{\mathcal A}$ (see Claim 5).
By Claim 4, each $X_{S}$, $S \in {\mathcal A}$, is a separable
unital projective subalgebra of $X$ and, by Claim 6, each limit
inclusion $i_{S} \colon X_{S} \to X$, $S \in {\mathcal A}$, is
doubly projective. This finishes the proof of the
implication (a) $\Longrightarrow$ (b).

Next we prove the implication (a) $\Longrightarrow$ (c). Since
$|T| = \tau$, we can write $T = \{ t_{\alpha} \colon \alpha < \tau\}$.
By Claim 5, for each $\alpha < \tau$ there exists a countable
admissible subset $S_{\alpha} \subseteq T$ such that
$t_{\alpha} \in S_{\alpha}$.
Let $T_{\alpha} = \cup\{ S_{\alpha} \colon \beta \leq \alpha\}$ and
$X_{\alpha} = X_{T_{\alpha}}$. Also let
$i_{\alpha}^{\alpha +1} \colon X_{\alpha} \to X_{\alpha +1}$ denote
the inclusion. Thus we have the well ordered continuous direct system
${\mathcal S}_{X} = \{ X_{\alpha}, i_{\alpha}^{\alpha +1}, \tau\}$.
It follows from the above constructions that
$X = \varinjlim{\mathcal S}_{X}$. According to Claims 3 and 4,
each $X_{\alpha}$, $\alpha < \tau$, is a unital projective
$C^{\ast}$-subalgebra of $X$. Since $T_{0} = S_{0}$ is countable,
we conclude that $X_{0}$ is separable. Claim 6 guarantees that
for each $\alpha < \tau$, both limit inclusions
$i_{\alpha} \colon X_{\alpha} \to X$ and
$i_{\alpha +1} \colon X_{\alpha +1} \to X$ are doubly
projective. Note that
$i_{\alpha} = i_{\alpha +1}\circ i_{\alpha}^{\alpha +1}$.
By Lemma \ref{L:inj}, $i_{\alpha +1}$ is a coretraction.
Consequently, by Lemma \ref{L:out}, $i_{\alpha}^{\alpha+1}$
is also doubly projective. Finally, in order to see that
$i_{\alpha}^{\alpha +1}$ has a separable type, note that
according to the above constructions and Lemma \ref{L:2},
we have the following commuting diagram

\[
\begin{CD}
X_{\alpha +1} @<\varphi_{T_{\alpha +1}}<<
\left(\bigstar_{{\mathbb C}}\{ X_{t} \colon t \in
T_{\alpha}\}\right)\bb \left(\bb\{ X_{t} \colon
t \in S_{\alpha +1}\}\right)\\
@A{i_{S}}AA @AA{\pi_{\bb\{ X_{t} \colon t \in T_{\alpha}\}}}A \\
X_{\alpha} @<\varphi_{T_{\alpha}}<<
\bigstar_{{\mathbb C}}\{ X_{t} \colon t \in T_{\alpha}\} ,
\end{CD}
\]

\noindent with surjective $\varphi_{T_{\alpha}}$ and
$\varphi_{T_{\alpha +1}}$. Clearly both
\[ \bigstar_{{\mathbb C}}\{ X_{t} \colon t \in T_{\alpha}\}\;\;
\text{and}\;\; \left(\bigstar_{{\mathbb C}}\{ X_{t} \colon t \in
T_{\alpha}\}\right)\bb \left(\bb\{ X_{t} \colon
t \in S_{\alpha +1}\}\right)\]

\noindent are projective
(as unital free products of projective $C^{\ast}$-algebras).
It only remains to note that since $S_{\alpha +1}$ is countable
and since each $X_{t}$ is separable, the unital free
product $\bb\{ X_{t} \colon
t \in S_{\alpha +1}\}$ is also separable
(\cite[Theorem 3.1]{blackadar1}). This completes the
proof of the implication (a) $\Longrightarrow$ (c).

In order to prove the implication (b) $\Longrightarrow$ (a)
observe that if
${\mathcal S}_{X} = \{ X_{\alpha}, i_{\alpha}^{\beta}, {\mathcal A}\}$
is a direct $C_{\omega}^{\ast}$-system satisfying properties
indicated in condition (b), then for any $\alpha \in {\mathcal A}$ the
$\alpha$-th limit inclusion $i_{\alpha} \colon X_{\alpha} \to X$ is doubly
projective and the $C^{\ast}$-algebra $X_{\alpha}$ is projective.
Consequently, by lemma \ref{L:iff}, $X$ is also projective.

Finally, the implication (c) $\Longrightarrow$ (a) follows
from Lemmas \ref{L:well} and \ref{L:iff}.
\end{proof}


\bigskip
\bigskip

\section{Basic properties of doubly projective
square diagrams and characterization of doubly projective
homomorphisms}\label{S:diagrams}

\subsection{Doubly projective diagrams}\label{SS:soublypd}

The pushout construction \cite{maclane} applied to the
category ${\mathcal C}_{1}^{\ast}$ leads us to the following definition
\cite{lope1}, \cite{pedersen1}. A commutative square diagram
$X_{1}X_{2}Y_{2}Y_{1}$,
consisting of unital $C^{\ast}$-algebras and unital
$\ast$-ho\-mo\-mor\-phisms,
is called {\em pushout}, if for any two coherent unital
$\ast$-homomorphisms
$g \colon X_{2} \to Z$ and $h \colon Y_{1} \to Z$ into
any unital $C^{\ast}$-algebra $Z$ (i.e. $g\circ i = h\circ f_{1}$),
there exists unique
unital $\ast$-homomorphism
$g \bigstar h \colon Y_{2} \to Z$ (a more informative notation $g\bigstar_{X_{1}}h$ for the sake of simplicity is replaced by $g \bigstar h$) such that
$(g \bigstar h)\circ f_{2} = g$ and $(g \bigstar h)\circ j = h$:

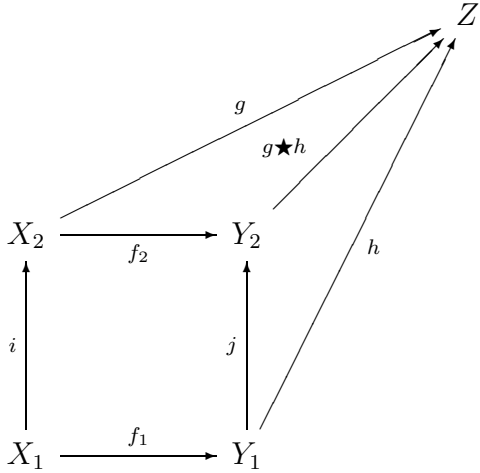
\[
\dgARROWLENGTH=2.0\dgARROWLENGTH
\begin{diagram}
\node[3]{Z}\\
\node{X_{2}} \arrow{e,b}{f_{2}} \arrow{ene,t}{g}
\node{Y_{2}} \arrow{ne,t,1}{g\bigstar h}\\
\node{X_{1}} \arrow{n,l}{i} \arrow{e,t}{f_{1}}
\node{Y_{1}} \arrow{n,l}{j} \arrow{nne,b}{h}
\end{diagram}
\]

\noindent $C^{\ast}$-algebra
$Y_{2}$ in such a case is isomorphic to the amalgamated free product $X_{2}\bigstar_{X_{1}}Y_{1}$, which is the quotient of the
unital free product
$X_{2}\bigstar_{{\mathbb C}}Y_{1}$ by the closed ideal
generated by $\{ i(x)-f_{1}(x) \colon x \in X_{1}\}$.
Also if $\pi \colon X_{2}\bigstar_{{\mathbb C}}Y_{1} \to Y_{2}$
denotes this quotient
homomorphism, then $\pi\circ \pi_{Y_{1}} = j$ and
$\pi \circ \pi_{X_{2}} = f_{2}$,
where $\pi_{Y_{1}} \colon Y_{1} \to X_{2}\bigstar_{{\mathcal C}}Y_{1}$ and
$\pi_{X_{2}} \colon X_{2} \to X_{2}\bigstar_{{\mathcal C}}Y_{1}$
denote canonical embeddings.
Here is the corresponding diagram

\[
\dgARROWLENGTH=2.0\dgARROWLENGTH
\begin{diagram}
\node[3]{X_{2}\bb Y_{1}} \arrow{sw,t,3}{\pi}\\
\node{X_{2}} \arrow{e,b}{f_{2}} \arrow{ene,t}{\pi_{X_{2}}}
\node{Y_{2}}\\
\node{X_{1}} \arrow{n,l}{i} \arrow{e,t}{f_{1}}
\node{Y_{1}} \arrow{n,l}{j} \arrow{nne,b}{\pi_{Y_{1}}}
\end{diagram}
\]

\begin{lem}\label{L:pushout}
Let 
\[
\begin{CD}
X_{2} @>f_{2}>> Y_{2}\\
@A{i}AA @AA{j}A\\
X_{1} @>f_{1}>> Y_{1}
\end{CD}
\]

\noindent be a pushout diagram,
consisting of unital $C^{\ast}$-algebras and
unital $\ast$-ho\-mo\-mor\-phisms. If $f_{1}$ is doubly projective,
then $f_{2}$ is also doubly projective.
\end{lem}
\begin{proof}
Let $p \colon A \to B$ be a surjective unital $\ast$-homomorphism
of unital $C^{\ast}$-algebras. Consider also two unital
$\ast$-homomorphisms
$g \colon X_{2} \to A$ and $h \colon Y_{2} \to B$ such
that $p\circ g = h\circ f_{2}$. Clearly 
\[ p\circ (g\circ i) =
(p \circ g)\circ i = (h\circ f_{2})\circ i =
(h\circ j)\circ f_{1} = h \circ (j \circ f_{1}) .\]
Since $f_{1}$ is doubly projective, there exists a unital
$\ast$-homomorphism
$\widetilde{k} \colon Y_{2} \to A$ such that
$g \circ i = \widetilde{k}\circ f_{1}$
and $h \circ j = p \circ \widetilde{k}$. Since the given diagram
is a pushout, we have a unital $\ast$-homomorphism
$k = g \bigstar \widetilde{k} \colon Y_{2} \to A$. Recall that
$g = (g \bigstar \widetilde{k}) \circ f_{2}$ and
$\widetilde{k} = (g \bigstar \widetilde{k})\circ j$.
Consequently it only remains to show that
$p \circ (g \bigstar \widetilde{k}) = h$. In order to
prove this equality
note that
\[ \left[ p \circ (g \bigstar \widetilde{k})\right] \circ f_{2} =
p \circ g = h \circ f_{2}\;\;\text{and}\;\;
\left[ p \circ (g \bigstar \widetilde{k})\right] \circ j =
p \circ \widetilde{k} = h \circ j .\]

Again, since the given diagram is a pushout, the above
equalities imply that
$p \circ (g \bigstar \widetilde{k}) = h$ as required.
\end{proof}

\begin{lem}\label{L:pushoutex}
Let $f \colon X \to Y$ be a unital $\ast$-homomorphism of unital
$C^{\ast}$-algebras. Let also $A$ be a unital $C^{\ast}$-algebra. Then
the diagram

\[
\begin{CD}
Y @>\pi_{Y}>> Y\bigstar_{{\mathbb C}}A\\
@A{f}AA @AA{(\pi_{Y}\circ f)\bigstar \operatorname{id}_{A}}A\\
X @>\pi_{X}>> X\bigstar_{{\mathbb C}}A
\end{CD}
\]
\noindent is a pushout.
\end{lem}
\begin{proof}
Consider the pushout

\[
\begin{CD}
Y @>\varphi>> Z\\
@A{f}AA @AA{\phi}A\\
X @>\pi_{X}>> X\bigstar_{{\mathbb C}}A
\end{CD}
\]

\noindent generated by the homomorphisms $f \colon X \to Y$ and
$\pi_{X} \colon X \to X\bigstar_{{\mathbb C}}A$. Since,
by the commutativity of the first diagram,

\begin{equation}\label{EQ:a}
\pi_{Y}\circ f = \left[ (\pi_{Y}\circ f)\bigstar
\operatorname{id}_{A}\right] \circ \pi_{X}\; ,
\end{equation}
it follows that there exists unique unital $\ast$-homomorphism
$p \colon Z \to Y\bigstar_{{\mathbb C}}A$ such that

\begin{equation}\label{EQ:aa}
p \circ \varphi = \pi_{Y}
\end{equation}

\noindent and

\begin{equation}\label{EQ:aaa}
p \circ \phi = (\pi_{Y}\circ f)\bigstar
\operatorname{id}_{A} .
\end{equation}

Let $\pi_{A} \colon A \to X\bigstar_{{\mathbb C}}A$ denote the canonical
injection of $A$ into the unital free product $X\bigstar_{{\mathbb C}}A$.
Consider the homomorphisms $\varphi \colon Y \to Z$ and
$\phi \circ \pi_{A} \colon A \to Z$. Since
$Y \bigstar_{{\mathbb C}}A$ is the unital free product,
there exists unique unital $\ast$-homomorphism
$q \colon Y \bigstar_{{\mathbb C}}A \to Z$ such that

\begin{equation}\label{EQ:b}
q \circ \pi_{Y} = \varphi
\end{equation}

\noindent and 

\begin{equation}\label{EQ:bb}
q \circ \lambda_{A} = \phi \circ \pi_{A}  ,
\end{equation}

\noindent where $\lambda_{A} \colon A \to Y \bigstar_{{\mathbb C}}A$
denotes the
canonical injection (not to be confused with $\pi_{A}$). Note that

\begin{equation}\label{EQ:bbb}
\lambda_{A} = \left[ (\pi_{Y}\circ f)\bigstar
\operatorname{id}_{A}\right] \circ \pi_{A} .
\end{equation}

In order to prove our statement we need to show
that $p$ is an
isomorphism. We accomplish this by proving that
$q\circ p = \operatorname{id}_{Z}$ and
$p \circ q = \operatorname{id}_{Y\bigstar_{{\mathbb C}}A}$.
The following diagram helps to visualize the situation.

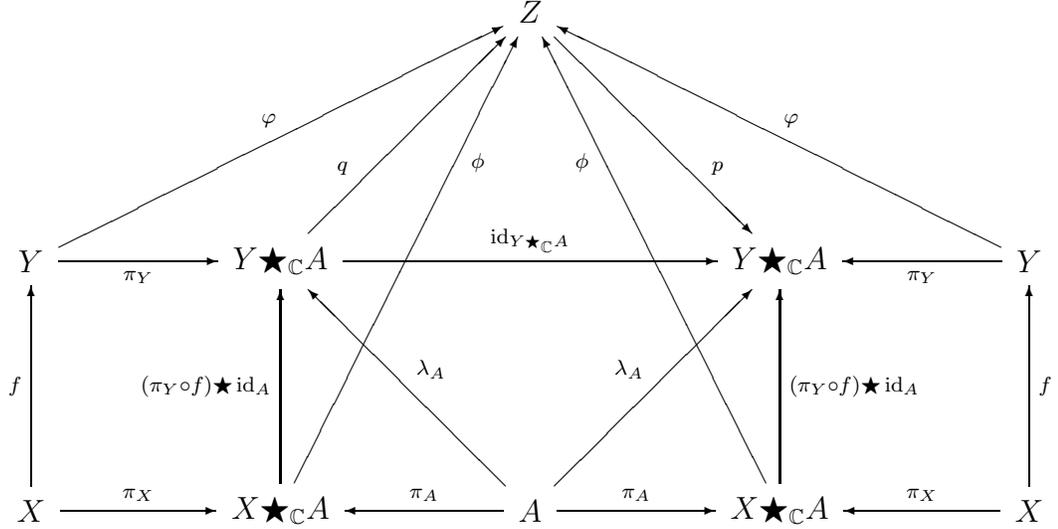
\[
\dgARROWLENGTH=2.0\dgARROWLENGTH
\begin{diagram}
\node[3]{Z} \arrow{se,t,3}{p}\\
\node{Y} \arrow{e,b}{\pi_{Y}} \arrow{ene,t}{\varphi}
\node{Y \bb A} \arrow{ne,t,1}{q} \arrow[2]{e,t}
{\operatorname{id}_{Y\bb A}}
\node[2]{Y\bb A} 
\node{Y} \arrow{w,b}{\pi_{Y}} \arrow{wnw,t}{\varphi}\\
\node{X} \arrow{n,l}{f} \arrow{e,t}{\pi_{X}}
\node{X\bb A} \arrow{n,l}{(\pi_{Y}\circ f)\bigstar
\operatorname{id}_{A}} \arrow{nne,b,3}{\phi}
\node{A} \arrow{w,t}{\pi_{A}} \arrow{nw,t}{\lambda_{A}} \arrow{e,t}{\pi_{A}} \arrow{ne,t}{\lambda_{A}}
\node{X \bb A} \arrow{n,r}{(\pi_{Y}\circ f)\bigstar
\operatorname{id}_{A}} \arrow{nnw,b,3}{\phi}
\node{X} \arrow{w,t}{\pi_{X}} \arrow{n,r}{f}
\end{diagram}
\]

First let us show that 

\begin{equation}\label{EQ:z}
q \circ \left[ (\pi_{Y}\circ f)\bigstar
\operatorname{id}_{A}\right]  = \phi .
\end{equation}
Since both $q \circ \left[ (\pi_{Y}\circ f)\bigstar
\operatorname{id}_{A}\right]$ and $\phi$ are defined on the
unital free product $X \bigstar_{{\mathbb C}}A$, (\ref{EQ:z}) will
be proved by examining compositions
of the above homomorphisms with $\pi_{A}$ and $\pi_{X}$. Observe that

\begin{equation}\label{EQ:zz}
\Bigl( q \circ \left[ (\pi_{Y}\circ f)\bigstar
\operatorname{id}_{A}\right]\Bigr) \circ \pi_{A}
\stackrel{(\ref{EQ:bbb})}{=} q\circ \lambda_{A}
\stackrel{(\ref{EQ:bb})}{=} \phi \circ \pi_{A}
\end{equation} 
\noindent and

\begin{equation}\label{EQ:zzz}
\Bigl( q \circ \left[ (\pi_{Y}\circ f)\bigstar
\operatorname{id}_{A}\right]\Bigr) \circ \pi_{X}
\stackrel{(\ref{EQ:a})}{=} q\circ \pi_{Y}\circ f
\stackrel{(\ref{EQ:b})}{=}\varphi \circ f = \phi \circ \pi_{X} .
\end{equation}

Note that (\ref{EQ:zz}) and (\ref{EQ:zzz}) imply (\ref{EQ:z}). 

Next note that

\begin{equation}\label{EQ:w}
q\circ p \circ \phi \stackrel{(\ref{EQ:aaa})}{=}
q \circ \left[ (\pi_{Y}\circ f)\bigstar
\operatorname{id}_{A}\right] \stackrel{(\ref{EQ:z})}{=} \phi 
\end{equation}

\noindent and

\begin{equation}\label{EQ:ww}
q \circ p \circ \varphi \stackrel{(\ref{EQ:aa})}{=}
q \circ \pi_{Y} \stackrel{(\ref{EQ:b})}{=} \varphi .
\end{equation}

Clearly (\ref{EQ:w}) and (\ref{EQ:ww}) imply the equality
$q\circ p = \operatorname{id}_{Z}$.

In order to establish the second equality
$p \circ q = \operatorname{id}_{Y\bigstar_{{\mathbb C}}A}$
we proceed in a similar way. Observe that

\begin{equation}\label{EQ:y}
p \circ q \circ \lambda_{A} \stackrel{(\ref{EQ:bb})}{=}
p \circ \phi \circ \pi_{A}= \left[ (\pi_{Y}\circ f)\bigstar
\operatorname{id}_{A}\right] \circ \pi_{A} = \lambda_{A}
\end{equation}

\noindent and

\begin{equation}\label{EQ:yy}
p \circ q \circ \pi_{Y} = p\circ \varphi = \pi_{Y} .
\end{equation}

As above, (\ref{EQ:y}) and (\ref{EQ:yy}) imply the required equality
$p \circ q = \operatorname{id}_{Y\bigstar_{{\mathbb C}}A}$.

This shows that $p$ is an isomorphism and completes the proof.
\end{proof}

\begin{lem}\label{L:seppushout}
Let 
\[
\begin{CD}
X_{2} @>f_{2}>> Y_{2}\\
@A{i}AA @AA{j}A\\
X_{1} @>f_{1}>> Y_{1}
\end{CD}
\]

\noindent be a pushout diagram,
consisting of projective unital $C^{\ast}$-algebras and
doubly projective ho\-mo\-mor\-phisms. If $f_{1}$ has a separable type,
then $f_{2}$ also has a separable type.
\end{lem}
\begin{proof}
Since $f_{1}$ has a separable type, we have the following commutative diagram

\[
\begin{CD}
Y_{1} @<\varphi_{Y_{1}}<<
X^{\prime}_{1}\bb Y^{\prime}_{1}\\
@A{f_{1}}AA @AA{\pi_{X^{\prime}_{1}}}A \\
X_{1} @<\varphi_{X_{1}}<<
X^{\prime}_{1},
\end{CD}
\]

\noindent where $X^{\prime}_{1}$ and $Y_{1}^{\prime}$ are
projective unital $C^{\ast}$-algebras, $Y_{1}^{\prime}$ in
addition is separable
and the unital $\ast$-homomorphisms $\varphi_{X_{1}}$ and
$\varphi_{Y_{1}}$ are surjective. By Lemma \ref{L:pushoutex},
the diagram

\[
\begin{CD}
X_{1}\bb Y_{1}^{\prime} @<(\pi_{X_{1}}\circ\varphi_{X_{1}})\bb \operatorname{id}_{Y_{1}^{\prime}}<<
X^{\prime}_{1}\bb Y^{\prime}_{1}\\
@A{\pi_{X_{1}}}AA @AA{\pi_{X^{\prime}_{1}}}A \\
X_{1} @<\varphi_{X_{1}}<<
X^{\prime}_{1},
\end{CD}
\]

\noindent is a pushout. Consequently there exists a unital
$\ast$-homomorphism $r \colon X_{1}\bb Y_{1}^{\prime} \to Y_{1}$
such that
$\varphi_{Y_{1}} = r\circ\left[ (\pi_{X_{1}}\circ\varphi_{X_{1}})\bb \operatorname{id}_{Y_{1}^{\prime}} \right]$. Since $\varphi_{Y_{1}}$
is surjective, the latter equality guarantees that $r$ is also surjective.
Thus we have the commutative diagram

\[
\begin{CD}
Y_{1} @<r<<
X_{1}\bb Y^{\prime}_{1}\\
@A{f_{1}}AA @AA{\pi_{X_{1}}}A \\
X_{1} @<\operatorname{id}_{X_{1}}<<
X_{1}.
\end{CD}
\]

Next consider the following diagram

\[
\dgARROWLENGTH=2.0\dgARROWLENGTH
\begin{diagram}
\node[3]{X_{2}\bb Y_{1}^{\prime}} \arrow{se,r,1}{\widetilde{r}}\\
\node{X_{2}} \arrow{ene,t}{\pi_{X_{2}}} \arrow{e,b}{f_{2}}
\node{Y_{2}} \arrow{ne,b,3}{\widetilde{s}} \arrow[2]{e,t,1}{\operatorname{id}_{Y_{2}}}
\node[2]{Y_{2}}
\node{X_{2}} \arrow{wnw,t}{\pi_{X_{2}}} \arrow{w,b}{f_{2}}\\
\node[3]{X_{1}\bb Y_{1}^{\prime}} \arrow{se,r}{r} \arrow[2]{n,r,1}{(\pi_{X_2}\circ i)\bb \operatorname{id}_{Y_{1}^{\prime}}}\\
\node{X_{1}} \arrow[2]{n,r}{i} \arrow{ene,t}{\pi_{X_{1}}} \arrow{e,b}{f_{1}}
\node{Y_{1}} \arrow[2]{n,l}{j} \arrow{ne,b}{s} \arrow[2]{e,t}{\operatorname{id}_{Y_{1}}}
\node[2]{Y_{1}} \arrow[2]{n,r}{j}
\node{X_{1}} \arrow[2]{n,l}{i} \arrow{wnw,t}{\pi_{X_{1}}} \arrow{w,b}{f_{1}}
\end{diagram}
\] 

\noindent in which, according to Lemma \ref{L:pushoutex}, the
subdiagram, represented by the back face of the above diagram,
is a pushout. Since $r$ is surjective and since $f_{1}$ is
doubly projective, there exists a unital $\ast$-homomorphism
$s \colon Y_{1} \to X_{1} \to Y_{1}^{\prime}$ such that
$r\circ s = \operatorname{id}_{Y_{1}}$ and $\pi_{X_{1}} = s\circ f_{1}$.
Now consider the unital $\ast$-homomorphisms
$j\circ r \colon X_{1}\bb Y_{1}^{\prime} \to Y_{2}$ and
$f_{2} \colon X_{2} \to Y_{2}$. Note that
$j\circ r \circ \pi_{X_{1}} = j\circ f_{1} = f_{2} \circ i$. Since,
as was indicated, the back face is a pushout, there exists the unique unital
$\ast$-homomorphism
$\widetilde{r} \colon X_{2}\bb Y_{1}^{\prime} \to Y_{2}$ such that
\[ \widetilde{r}\circ\left[ (\pi_{X_{2}}\circ i)\bb
\operatorname{id}_{Y_{1}^{\prime}}\right] = j\circ r\;
\text{and}\; \widetilde{r}\circ \pi_{X_{2}} = f_{2} .\]

It only remains to show that $\widetilde{r}$ is surjective. To see this
consider the homomorphisms $\left[ (\pi_{X_{2}}\circ i)\bb \operatorname{id}_{Y_{1}^{\prime}}\right]\circ s \colon Y_{1}
\to X_{2} \bb Y_{1}^{\prime}$ and $\pi_{X_{2}} \colon X_{2}
\to X_{2}\bb Y_{1}^{\prime}$. Clearly
\[ \left[ (\pi_{X_{2}}\circ i)\bb
\operatorname{id}_{Y_{1}^{\prime}}\right]\circ s \circ f_{1}= 
\left[ (\pi_{X_{2}}\circ i)\bb
\operatorname{id}_{Y_{1}^{\prime}}\right]\circ \pi_{X_{1}} =
\pi_{X_{2}}\circ i .\]

\noindent Since the originally given diagram is a pushout, there
exists a unital $\ast$-ho\-mo\-mor\-phism
$\widetilde{s} \colon Y_{2} \to X_{2}\bb Y_{1}^{\prime}$ such that 
$\widetilde{s}\circ j = \left[ (\pi_{X_{2}}\circ i)\bb
\operatorname{id}_{Y_{1}^{\prime}}\right]\circ s$ and
$\widetilde{s}\circ f_{2} = \pi_{X_{2}}$. Straightforward
verification (based on the universality properties of the two
pushout diagrams involved) shows that $\widetilde{r}\circ\widetilde{s} = \operatorname{id}_{Y_{2}}$. This suffices to conclude that
$\widetilde{r}$ is surjective. Consequently the
homomorphism $f_{2}$ has a separable type.
\end{proof}

\begin{defin}\label{D:characteristic}
A characteristic $\ast$-homomorphism of a commutative square
diagram $X_{1}X_{2}Y_{1}Y_{2}$ is the
$\ast$-homomorphism $\chi = f_{2}\bigstar j$

\bigskip

\[
\dgARROWLENGTH=2.0\dgARROWLENGTH
\begin{diagram}
\node[3]{X_{2}\bigstar_{X_{1}} Y_{1}} \arrow{sw,t,3}{\chi}\\
\node{X_{2}} \arrow{e,b}{f_{2}} \arrow{ene,t}{\varphi_{X_{2}}}
\node{Y_{2}}\\
\node{X_{1}} \arrow{n,l}{i} \arrow{e,t}{f_{1}}
\node{Y_{1}} \arrow{n,l}{j} \arrow{nne,b}{\varphi_{Y_{1}}}
\end{diagram}
\]
\end{defin}

Note that a commutative square diagram is a pushout if and only
if its characteristic $\ast$-homomorphism is an isomorphism.

\begin{defin}\label{D:doublydiagr}
A commutative square diagram, consisting of
unital $C^{\ast}$-algebras
and unital $\ast$-homomorphisms, is called
doubly projective, if its
characteristic $\ast$-homomorphism is doubly projective.
\end{defin}

\begin{lem}\label{L:doublydiagr}
Let 

\[
\begin{CD}
X_{2} @>f_{2}>> Y_{2}\\
@A{i}AA @AA{j}A\\
X_{1} @>f_{1}>> Y_{1}
\end{CD}
\]
\noindent be a doubly projective square diagram.
If $f_{1}$ is doubly projective, then $f_{2}$ is also doubly
projective. Moreover,
for any unital surjective $\ast$-homomorphism $p \colon A \to B$ of
unital $C^{\ast}$-algebras and any three unital
$\ast$-ho\-mo\-mor\-phisms
$g \colon X_{2} \to A$, $h \colon Y_{2} \to B$ and
$k_{1} \colon Y_{1} \to A$ such that $h\circ f_{2} = p\circ g$,
$g\circ i = k_{1}\circ f_{1}$ and $h\circ j = p\circ k_{1}$,

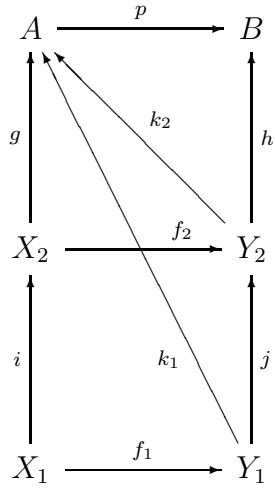
\[
\dgARROWLENGTH=2.0\dgARROWLENGTH
\begin{diagram}
\node{A} \arrow{e,t}{p} \node{B}\\
\node{X_{2}} \arrow{n,l}{g} \arrow{e,t,3}{f_{2}}
\node{Y_{2}} \arrow{n,r}{h} \arrow{nw,t}{k_{2}}\\
\node{X_{1}} \arrow{n,l}{i} \arrow{e,t}{f_{1}}
\node{Y_{1}} \arrow{n,r}{j} \arrow{nnw,b,1}{k_{1}}
\end{diagram} 
\]

\noindent there exists a unital $\ast$-homomorphism
$k_{2} \colon Y_{2} \to A$
such that $k_{2}\circ f_{2} = g$, $p\circ k_{2} = h$ and
$k_{2}\circ j = k_{1}$.
\end{lem}
\begin{proof}
Consider the pushout diagram
$X_{1}X_{2}X_{2}\bigstar_{X_{1}}Y_{2}Y_{2}$ generated by
the $\ast$-homomorphisms 
$i \colon X_{1} \to X_{2}$ and $f_{1} \colon X_{1} \to Y_{1}$.
Since $f_{1}$ is doubly projective, it follows, by Lemma
\ref{L:pushout}, that $\varphi_{X_{2}}$ is also doubly projective.
Since the characteristic $\ast$-homomorphism
$\chi \colon X_{2}\bigstar_{X_{1}}X_{1} \to Y_{2}$ of the
originally given diagram is doubly projective, it follows,
by Lemma \ref{L:comp},
that the composition $f_{2} = \chi \circ \varphi_{X_{2}}$
is doubly projective.
This proves the first part of our statement.

In order to prove the second part of Lemma consider the
following diagram in which all objects satisfy the above
formulated assumptions:\\

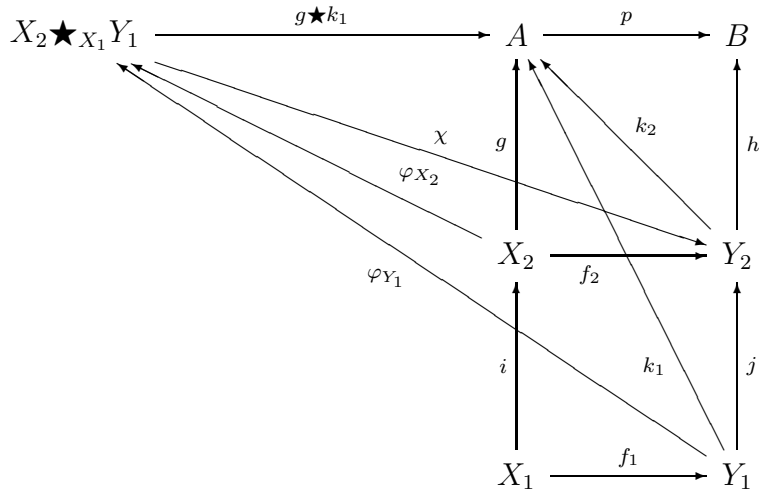
\[
\dgARROWLENGTH=2.0\dgARROWLENGTH
\begin{diagram}
\node{X_{2}\bigstar_{X_{1}}Y_{1}} \arrow{kk,t}{\chi} \arrow[2]{e,t}{g\bigstar k_{1}}  \node[2]{A}  \arrow{e,t}{p} \node{B}\\
\node[3]{X_{2}} \arrow{wnw,l,1}{\varphi_{X_{2}}} \arrow{n,l}{g} \arrow{e,b,1}{f_{2}}
\node{Y_{2}}  \arrow{n,r}{h} \arrow{nw,t}{k_{2}}\\
\node[3]{X_{1}} \arrow{n,l}{i} \arrow{e,t}{f_{1}}
\node{Y_{1}} \arrow{k,b}{\varphi_{Y_{1}}} \arrow{n,r}{j} \arrow{nnw,b,1}{k_{1}}
\end{diagram} 
\]

Since $X_{1}X_{2}X_{2}\bigstar_{X_{1}}Y_{2}Y_{2}$ is a
pushout diagram
and since the $\ast$-homomorphisms
$g \colon X_{2} \to A$ and $k_{1} \colon Y_{1} \to A$
satisfy the equality
$g \circ i = k_{1}\circ f_{1}$, there exists unique
$\ast$-homomorphism
$g\bigstar k_{1} \colon X_{2}\bigstar_{X_{1}}Y_{1} \to A$
such that 
\begin{equation}\label{EQ:t}
g = (g\bigstar k_{1}) \circ \varphi_{X_{2}}
\end{equation}

\noindent and
\begin{equation}\label{EQ:tt}
k_{1} = (g\bigstar k_{1}) \circ \varphi_{Y_{1}} .
\end{equation}

In order to prove that $p \circ (g\bigstar k_{1}) = h\circ \chi$,
first observe that
\begin{equation}\label{EQ:pp}
 \left[ p \circ (g\bigstar k_{1})\right]\circ \varphi_{Y_{1}} =
p\circ k_{1} = h\circ j = \left[ h\circ \chi\right] \circ \varphi_{Y_{1}} .
\end{equation}
\noindent Secondly,
\begin{equation}\label{EQ:ppp}
 \left[ p \circ (g\bigstar k_{1})\right] \circ \varphi_{X_{2}} = p\circ g =
h\circ f_{2} = \left[ h\circ \chi\right] \circ\varphi_{X_{2}} .
\end{equation}
Since $X_{1}X_{2}X_{2}\bigstar_{X_{1}}Y_{2}Y_{2}$ is a
pushout diagram, (\ref{EQ:pp}) and (\ref{EQ:ppp}), imply the required equality
$p \circ (g\bigstar k_{1}) = h\circ \chi$.

Since $\chi$ is doubly projective the latter equality guarantees
the existence of a unital $\ast$-homomorphism $k_{2} \colon Y_{2}\to A$
such that $p\circ k_{2} = h$ and $k_{2}\circ \chi = g\bigstar k_{1}$. 
The straitforward verification shows that $k_{2}\circ f_{2} = g$ and
$k_{2}\circ j = k_{1}$

\[ k_{2}\circ f_{2} = k_{2}\circ \chi \circ
\varphi_{X_{2}} = (g\bigstar k_{1})\circ
\varphi_{X_{2}}\stackrel{(\ref{EQ:t})}{=}g \]
\noindent and
\[ k_{2}\circ j = k_{2}\circ \chi \circ \varphi_{Y_{1}} = (g\bigstar k_{1})\circ
\varphi_{Y_{1}}\stackrel{(\ref{EQ:tt})}{=}k_{1} .\]
This completes the proof of Lemma \ref{L:doublydiagr}.
\end{proof}

\begin{pro}\label{P:doubly}
Let ${\mathcal S}_{X} = \{ X_{\alpha}, i_{\alpha}^{\alpha +1}, \tau\}$ and
${\mathcal S}_{Y} = \{ Y_{\alpha}, j_{\alpha}^{\alpha +1}, \tau\}$ be two
well ordered continuous direct systems consisting of unital $C^{\ast}$-algebras
and unital $\ast$-homomorphisms. Let
\[ \{ f_{\alpha} \colon X_{\alpha} \to Y_{\alpha};
\alpha \in \tau\} \colon {\mathcal S}_{X} \to {\mathcal S}_{Y}\]

\noindent be a morphism between these systems such that all arising adjacent
square diagrams 

\[
\begin{CD}
X_{\alpha +1} @>f_{\alpha +1}>> Y_{\alpha +1}\\
@A{i_{\alpha}^{\alpha +1}}AA @AA{j_{\alpha}^{\alpha +1}}A\\
X_{\alpha} @>f_{\alpha}>> Y_{\alpha}
\end{CD}
\]

\noindent are doubly projective. If $f_{0} \colon X_{0} \to Y_{0}$
is doubly projective, then
the limit homomorphism
$\varinjlim\{ f_{\alpha}\} \colon \varinjlim{\mathcal S}_{X}
\to {\mathcal S}_{Y}$ is also doubly projective.
\end{pro}
\begin{proof}
Let $p \colon A \to B$ be a unital surjective $\ast$-homomorphism
of unital $C^{\ast}$-algebras. Consider two unital
$\ast$-homomorphisms
\[ g \colon \varinjlim{\mathcal S}_{X} \to A\;\;\text{and}\;\;
h \colon \varinjlim{\mathcal S}_{Y} \to B \]

\noindent such that $p \circ g = h \circ \varinjlim\{ f_{\alpha}\}$. 
Our goal is to construct a unital $\ast$-homomorphism 
\[ k \colon \varinjlim{\mathcal S}_{Y} \to A \]

\noindent such that $k \circ \varinjlim\{ f_{\alpha}\} = g$ and
$p \circ k = h$.
Let 
\[ g_{\alpha} = g\circ i_{\alpha} \colon X_{\alpha} \to A \;\;
\text{and}\;\; h_{\alpha} = h\circ j_{\alpha} \colon Y_{\alpha} \to B ,
\alpha < \tau .\]

We now construct (by induction) a collection of unital $\ast$-homomorphisms
\[ k_{\alpha} \colon  Y_{\alpha} \to A ,\; \alpha < \tau ,\]
so that the following conditions are satisfied:
\begin{itemize}
\item[(a)]
$g_{\alpha} = k_{\alpha}\circ f_{\alpha}$, $\alpha < \tau$.
\item[(b)]
$h_{\alpha} = p \circ k_{\alpha} $, $\alpha < \tau$.
\item[(c)]
$k_{\alpha} = k_{\alpha +1}\circ j_{\alpha}^{\alpha +1}$, $\alpha < \tau$.
\item[(d)]
$k_{\alpha} = \varinjlim\{ k_{\beta}; \beta < \alpha\}$, whenever $\alpha$ is a
limit ordinal number with $\alpha < \tau$.
\end{itemize}

By our assumption, the $\ast$-homomorphism $f_{0}$ is doubly
projective. Consequently there exists a unital $\ast$-homomorphism
$k_{0} \colon Y_{0} \to A$ such that $g_{0} = k_{0}\circ f_{0}$ and
$h_{0} = p \circ k_{0}$.

Suppose that for each $\alpha < \gamma$, where $\gamma < \tau$, we
have already constructed unital $\ast$-homomorphisms
$k_{\alpha} \colon Y_{\alpha} \to A$ satisfying conditions
(a)--(d) for appropriate indices. Let us construct a unital $\ast$-homomorphism
$k_{\gamma} \colon Y_{\gamma} \to A$.

If $\gamma$ is a limit ordinal number, then let (consult
with Subsection \ref{SS:morphisms})

\[ k_{\gamma} = \varinjlim\{ f_{\alpha} \colon \alpha < \gamma \} .\]

The continuity of the direct systems ${\mathcal S}_{X}$ and
${\mathcal S}_{Y}$ guarantees that
$g_{\gamma} = k_{\gamma}\circ f_{\gamma}$, $h_{\gamma} = p \circ k_{\gamma}$
and $k_{\alpha} = j_{\alpha}^{\gamma}\circ k_{\gamma}$ for each $\alpha < \gamma$.

If $\gamma = \alpha +1$, then, by the assumption, the diagram

\[
\begin{CD}
X_{\alpha +1} @>f_{\alpha +1}>> Y_{\alpha +1}\\
@A{i_{\alpha}^{\alpha +1}}AA @AA{j_{\alpha}^{\alpha +1}}A\\
X_{\alpha} @>f_{\alpha}>> Y_{\alpha}
\end{CD}
\]

\noindent is doubly projective. Therefore, by Lemma
\ref{L:doublydiagr}, there exists a unital $\ast$-ho\-mo\-mor\-phism
$k_{\alpha +1} \colon Y_{\alpha +1} \to A$ such that
$g_{\alpha +1} = k_{\alpha +1}\circ f_{\alpha +1}$,
$h_{\alpha +1} = p \circ k_{\alpha +1}$
and $k_{\alpha} = j_{\alpha}^{\alpha +1}\circ k_{\alpha +1}$.

Thus, the $\ast$-homomorphisms $k_{\alpha}$ are now constructed
for each $\alpha < \tau$. It only remains to note that the
$\ast$-homomorphism
$k = \varinjlim\{ k_{\alpha}\} \colon \varinjlim{\mathcal S}_{Y} \to A$
satisfies all the required properties.
\end{proof}


\subsection{Characterization of doubly projective
homomorphisms}\label{SS:charact}
\begin{thm}\label{T:q}
Let $f \colon X \to Y$ be a unital $\ast$-homomorphism
between unital $C^{\ast}$-algebras of the same
density. Then $f$ is doubly projective homomorphism of separable type if
and only if there exist direct $C_{\omega}^{\ast}$-systems
${\mathcal S}_{X} = \{ X_{\alpha}, i_{\alpha}, {\mathcal A}\}$,
${\mathcal S}_{Y} = \{ Y_{\alpha}, j_{\alpha}^{\beta}, {\mathcal A}\}$
and a morphism $\{ f_{\alpha} \colon X_{\alpha} \to
Y_{\alpha}; \alpha \in {\mathcal A}\} \colon {\mathcal S}_{X}
\to {\mathcal S}_{Y}$, satisfying the following conditions:
\begin{itemize}
\item[(a)]
The indexing set ${\mathcal A}$ is cofinal and $\omega$-closed
in $\exp_{\omega}\tau$.
\item[(b)]
$X = \varinjlim{\mathcal S}_{X}$,
$Y = \varinjlim{\mathcal S}_{Y}$,
$f = \varinjlim\{ f_{\alpha}; \alpha \in {\mathcal A}\}$.
\item[(c)]
$X_{\alpha}$ and $Y_{\alpha}$ are separable unital projective
$C^{\ast}$-algebras, $\alpha \in {\mathcal A}$.
\item[(d)]
The $\alpha$-th limit inclusions
$i_{\alpha} \colon X_{\alpha} \to X$ and
$j_{\alpha} \colon Y_{\alpha} \to Y$ are doubly
projective, $\alpha \in {\mathcal A}$.
\item[(e)]
$f_{\alpha} \colon X_{\alpha} \to Y_{\alpha}$ is doubly
projective, $\alpha \in {\mathcal A}$.
\item[(f)]
All $\alpha$-th limit diagrams ($\alpha \in {\mathcal A}$)
\[
\begin{CD}
X @>f>> Y\\
@A{i_{\alpha}}AA @AA{j_{\alpha}}A\\
X_{\alpha} @>f_{0}>> Y_{\alpha},
\end{CD}
\]

\noindent are pushouts.
\end{itemize}
\end{thm}
\begin{proof}
{\em Part I}. Let $f \colon X \to Y$ be a doubly
projective homomorphism of separable type. We will show
the existence of the above indicated direct
$C_{\omega}^{\ast}$-systems and of a morphism, satisfying
the required properties.

If $Y$ is separable, then the statement is trivial. Indeed, by
Lemma \ref{L:inj}, $f$ is injective and consequently $X$ is also
separable. Let $X_{0} = X$, $Y_{0} = Y$, $p = \operatorname{id}_{X}$,
$q = \operatorname{id}_{Y}$ and $f_{0} = f$. Obviously the diagram

\[
\begin{CD}
X @>f>> Y\\
@A{\operatorname{id}_{X}}AA @AA{\operatorname{id}_{Y}}A\\
X @>f>> Y,
\end{CD}
\]

\noindent is a pushout.

Now consider the case $d(Y) = \tau > \omega$. By our assumption,
the homomorphism $f$ has a separable type. This means
(see Definition \ref{D:septype}) that there exist a
projective unital $C^{\ast}$-algebra $Z$ such that $d(Z) = d(X)$,
a separable
projective unital $C^{\ast}$-algebra $K$ and two
surjective unital $\ast$-homomorphisms
$\varphi_{X} \colon Z \to X$ and
$\varphi_{Y} \colon Z \bb K \to Y$ such that
$f \circ \varphi_{X} = \varphi_{Y}\circ \pi_{Z}$, where
$\pi_{Z} \colon Z \to  Z \bb K$
denotes the natural inclusion. In other words, the following diagram

\[
\begin{CD}
Y @<\varphi_{Y}<<
Z\bb K\\
@A{f}AA @AA{\pi_{Z}}A \\
X @<\varphi_{X}<< Z
\end{CD}
\]

\noindent commutes. 

Since $X$ is projective and $\varphi_{X} \colon Z \to X$ is
surjective, there exists a unital $\ast$-homomorphism
$\phi_{X} \colon X \to Z$ such that
$\varphi_{X}\circ\phi_{X} = \operatorname{id}_{X}$.
Now consider the
square diagram

\[
\dgARROWLENGTH=2.0\dgARROWLENGTH
\begin{diagram}
\node{Y} \node{Y} \arrow{w,t}{\operatorname{id}_{Y}} \arrow{sw,t}{\phi_{Y}}\\
\node{Z\bb K} \arrow{n,l}{\varphi_{Y}} \node{X}
\arrow{w,t}{\pi_{Z}\circ \phi_{X}} \arrow{n,r}{f} 
\end{diagram}
\]

\noindent which obviously commutes. To see this note that
\[ \varphi_{Y}\circ \pi_{Z}\circ \phi_{X} =
f\circ \varphi_{X}\circ\phi_{X} = f .\]

Since $\varphi_{Y}$ is surjective and since $f$ is doubly projective,
there exists a unital $\ast$-homomorphism
$\phi_{Y} \colon Z \bb K \to Y$ (indicated in the above diagram
as the diagonal arrow) such that
$\varphi_{Y}\circ \phi_{Y} = \operatorname{id}_{Y}$ and
$\phi_{Y}\circ f = \pi_{Z}\circ \phi_{X}$. Thus we have the
commutative diagram

\[
\begin{CD}
Y @>\phi_{Y}>> Z\bb K\\
@A{f}AA @AA{\pi_{Z}}A \\
X @>\phi_{X}>> Z.
\end{CD}
\]

Next observe that the $C^{\ast}$-algebras $X$, $Y$,
$Z$ and $Z \bb K$ all have density $\leq \tau$. Consequently, by
Theorem \ref{T:charact}, $X = \varinjlim{\mathcal S}_{X}$,
$Y = \varinjlim{\mathcal S}_{Y}$ and $Z = \varinjlim{\mathcal S}_{Z}$, where
${\mathcal S}_{X} = \{ X_{\alpha}, i_{\alpha}^{\beta}, {\mathcal A}_{X}\}$,
${\mathcal S}_{Y} = \{ Y_{\alpha}, j_{\alpha}^{\beta}, {\mathcal A}_{Y}\}$
and ${\mathcal S}_{Z} = \{ Z_{\alpha}, s_{\alpha}^{\beta}, {\mathcal A}_{Z}\}$
are direct $C^{\ast}_{\omega}$-systems 
consisting of separable unital projective $C^{\ast}$-algebras and
doubly projective
limit inclusions
$i_{\alpha} \colon X_{\alpha} \to X$, $\alpha \in {\mathcal A}_{X}$,
$j_{\alpha} \colon Y_{\alpha} \to Y$, $\alpha \in {\mathcal A}_{Y}$,
and $s_{\alpha} \colon Z_{\alpha} \to Z$, $\alpha \in {\mathcal A}_{Z}$.
Also note that all three indexing sets ${\mathcal A}_{X}$,
${\mathcal A}_{Y}$
and ${\mathcal A}_{Z}$ are cofinal and $\omega$-closed subsets
of $\exp_{\omega}\tau$. Next observe that the unital free
product $Z \bb K$
is also the limit of the direct system
${\mathcal S}_{Z\bb K} = \{ Z_{\alpha}\bb K,
s_{\alpha}^{\beta}\bb \operatorname{id}_{K}, {\mathcal A}_{Z}\}$
(straightforward 
verification using the universality properties of unital free
products and limits of
direct systems; see also Section \ref{S:free}). An important
consequence
of the fact that $f$ has a separable type is that $K$ is a
separable $C^{\ast}$-algebra. This guarantees, according to
\cite[Theorem 3.1]{blackadar1}, that each $C^{\ast}$-algebra
$Z_{\alpha}\bb K$, $\alpha \in {\mathcal A}_{Z}$,
is separable and, as a result, ${\mathcal S}_{Z\bb K}$ is actually
a direct $C_{\omega}^{\ast}$-system.

For each $\alpha \in {\mathcal A}_{Z}$ let $\widetilde{X}_{\alpha} = \operatorname{cl}_{X}\left(\varphi_{X}(Z_{\alpha})\right)$. Let also $\widetilde{i}_{\alpha}^{\beta} \colon \widetilde{X}_{\alpha}
\to \widetilde{X}_{\beta}$, $\alpha \leq \beta$,
$\alpha ,\beta \in {\mathcal A}_{Z}$ denote the corresponding
inclusion. Similarly, for each
$\alpha \in {\mathcal A}_{Z}$ let $\widetilde{Y}_{\alpha} = \operatorname{cl}_{Y}\left(\varphi_{Y}(Z_{\alpha}\bb K)\right)$ and 
$\widetilde{j}_{\alpha}^{\beta} \colon
\widetilde{Y}_{\alpha} \to \widetilde{Y}_{\beta}$, $\alpha \leq \beta$,
$\alpha ,\beta \in {\mathcal A}_{Z}$ denote the corresponding
inclusion. It is easy to see that the systems
$\widetilde{\mathcal S}_{X} =
\{ \widetilde{X}_{\alpha}, \widetilde{i}_{\alpha}^{\beta},
{\mathcal A}_{Z}\}$
and $\widetilde{\mathcal S}_{Y} =
\{ \widetilde{Y}_{\alpha}, \widetilde{j}_{\alpha}^{\beta},
{\mathcal A}_{Z}\}$
are direct $C_{\omega}^{\ast}$-systems such that
$\varinjlim\widetilde{\mathcal S}_{X} = X$ and
$\varinjlim\widetilde{\mathcal S}_{Y} = Y$.

Since the indexing sets ${\mathcal A}_{X}$, ${\mathcal A}_{Y}$
and ${\mathcal A}_{Z}$ are cofinal and $\omega$-closed in
$\exp_{\omega}\tau$, we can conclude, by Proposition \ref{P:3.1.1},
that the intersection
${\mathcal B} = {\mathcal A}_{X} \cap {\mathcal A}_{Y}
\cap {\mathcal A}_{Z}$
is still cofinal and $\omega$-closed in $\exp_{\omega}\tau$.

Next we consider six homomorphisms 
\[ \varphi_{X} \colon \varinjlim{\mathcal S}_{Z}|{\mathcal B}
\to \varinjlim{\mathcal S}_{X}|{\mathcal B}, \;\;\;\;\;\;\;\;\;\;\;\;\;\;\;\;\;\;\;\;\;\;\; \varphi_{Y}
\colon \varinjlim{\mathcal S}_{Z\bb K}|{\mathcal B} \to
\varinjlim{\mathcal S}_{Y}|{\mathcal B},\]
\[ \phi_{X} \colon \varinjlim{\mathcal S}_{X}|{\mathcal B}
\to \varinjlim{\mathcal S}_{Z}|{\mathcal B},\;\;\;\;\;\;\;\;\;\;\;\;\;\;\;\;\;\;\;\;\;\;\;
\phi_{Y} \colon \varinjlim{\mathcal S}_{Y}|{\mathcal B}
\to \varinjlim{\mathcal S}_{Z\bb K}|{\mathcal B},
\]
\[ \pi_{Z} \colon \varinjlim{\mathcal S}_{Z}|{\mathcal B}
\to \varinjlim{\mathcal S}_{Z\bb K}|{\mathcal B}
\;\;\;\;\;\;\;\;\;\;\text{and}\;\;\;\;\;\;\;\;\;\; f
\colon \varinjlim{\mathcal S}_{X}|
{\mathcal B} \to \varinjlim{\mathcal S}_{Y}|{\mathcal B} .\]

Three of these homomorphisms are, by construction, the limits
of associated morphisms

\[ \varphi_{X} = \varinjlim\{ \varphi_{X}^{\alpha} \colon
Z_{\alpha} \to X_{\alpha}; {\mathcal B}\} ,\;\;\text{where}\;\;
\varphi_{X}^{\alpha} = \varphi_{X}|Z_{\alpha} , \;\;
\alpha \in {\mathcal B} ,\]

\[ \varphi_{Y} = \varinjlim\{ \varphi_{Y}^{\alpha} \colon
Z_{\alpha}\bb K \to Y_{\alpha}; {\mathcal B}\} ,\;\;
\text{where}\;\; \varphi_{Y}^{\alpha} = \varphi_{Y}|
\left(Z_{\alpha}\bb K\right) ,\;\;\alpha \in {\mathcal B} ,\]

\noindent and

\[ \pi_{Z} = \varinjlim\{ \pi_{Z_{\alpha}} \colon
Z_{\alpha} \to Z_{\alpha} \bb K; {\mathcal B}\} ,\;
\text{where}\; \pi_{Z_{\alpha}}\;\text{is the canonical
inclusion},\;\alpha \in {\mathcal B} ,\]

We apply Theorem \ref{T:spectral} to the remaining three
homomorphisms $\varphi_{X}$, $\varphi_{Y}$ and $f$ and
conclude that there exist cofinal and $\omega$-complete subsets
${\mathcal B}_{\varphi_{X}}$, ${\mathcal B}_{\varphi_{Y}}$ and
${\mathcal B}_{f}$ of ${\mathcal B}$ and morphisms
\[ \{ \varphi_{X}^{\alpha} \colon X_{\alpha} \to Z_{\alpha};
{\mathcal B}_{\varphi_{X}}\} \colon {\mathcal S}_{X}|
{\mathcal B}_{\varphi_{X}} \to {\mathcal S}_{Z}|
{\mathcal B}_{\varphi_{X}} ,\]
\[ \{ \varphi_{Y}^{\alpha} \colon Y_{\alpha} \to Z_{\alpha}\bb K;
{\mathcal B}_{\varphi_{Y}}\} \colon {\mathcal S}_{Y}|
{\mathcal B}_{\varphi_{X}} \to {\mathcal S}_{Z\bb K}|
{\mathcal B}_{\varphi_{Y}} ,\]
\noindent and
\[ \{ f_{\alpha} \colon X_{\alpha} \to Y_{\alpha};
{\mathcal B}_{f}\} \colon {\mathcal S}_{X}|
{\mathcal B}_{f} \to {\mathcal S}_{Y}|
{\mathcal B}_{f} \]

\noindent such that
\[ \varphi_{X} = \varinjlim\{ \varphi_{X}^{\alpha}; \alpha
\in {\mathcal B}_{\varphi_{X}}\} ,\;
\varphi_{Y} = \varinjlim\{ \varphi_{Y}^{\alpha};
\alpha \in {\mathcal B}_{\varphi_{Y}}\}\;\text{and} \;
f = \varinjlim\{ f_{\alpha} ; \alpha \in {\mathcal B}_{f}\} .\]

Note that, by Proposition \ref{P:3.1.1}, the intersection
${\mathcal A} = {\mathcal B}_{\varphi_{X}} \cap
{\mathcal B}_{\varphi_{Y}} \cap {\mathcal B}_{f}$ is cofinal
and $\omega$-closed in ${\mathcal B}$ (and consequently
in $\exp_{\omega}\tau$).

For each $\alpha \in {\mathcal A}$ we have the
following commutative diagram:

\[
\dgARROWLENGTH=1.2\dgARROWLENGTH
\begin{diagram}
\node{X} \node[2]{Z}  \arrow[2]{w,t}{\varphi_{X}}\arrow[2]{e,t}{\pi_{Z}} \node[2]{Z \bb K} \arrow{e,t}{\varphi_{Y}} \node{Y}\\
\node[2]{X} \arrow{nw,l,<>}{\operatorname{id}_{X}}\arrow{ne,l}{\phi_{X}} \arrow[2]{e,t,3}{f} \node[2]{Y} \arrow{ne,l}{\phi_{Y}} \arrow{ene,r,<>,3}{\operatorname{id}_{Y}}\\
\node{X_{\alpha}} \arrow[2]{n,l}{i_{\alpha}} \node[2]{Z_{\alpha}} \arrow[2]{w,t,1}{\varphi_{X_{\alpha}}}
\arrow[2]{n,r,3}{s_{\alpha}} \arrow[2]{e,t,1}{\pi_{Z_{\alpha}}} \node[2]{Z_{\alpha} \bb K} \arrow{e,t}{\varphi_{Y}^{\alpha}} \arrow[2]{n,r}{s_{\alpha}\bb \operatorname{id}_{K}} \node{Y_{\alpha}}
\arrow[2]{n,r}{j_{\alpha}} \\
\node[2]{X_{\alpha}} \arrow{nw,r,<>}{\operatorname{id}_{X_{\alpha}}} \arrow[2]{n,l,3}{i_{\alpha}} \arrow{ne,l}{\phi_{X}^{\alpha}}
\arrow[2]{e,t}{f_{\alpha}} \node[2]{Y_{\alpha}} \arrow[2]{n,r,3}{j_{\alpha}}
\arrow{ne,l}{\phi_{Y}^{\alpha}} \arrow{ene,r,<>}{\operatorname{id}_{Y_{\alpha}}}
\end{diagram}
\]

\noindent Note that, by Theorem \ref{T:charact}, we may
without loss of generality assume that the limit inclusions
$i_{\alpha} \colon X_{\alpha} \to X$ and
$j_{\alpha} \colon Y_{\alpha} \to Y$, $\alpha \in {\mathcal A}$,
are doubly projective. This observation coupled with Lemma \ref{L:out}
guarantees that the homomorphism
$f_{\alpha} \colon X_{\alpha} \to Y_{\alpha}$,
$\alpha \in {\mathcal A}$, is also doubly projective.

It is now clear that in order to complete the proof it suffices
to show that the diagram (the front face of the above cubic diagram)

\[
\begin{CD}
X @>f>> Y\\
@A{i_{\alpha}}AA @AA{j_{\alpha}}A \\
X_{\alpha} @>f{\alpha}>> Y_{\alpha}
\end{CD}
\]

\noindent is a pushout for an arbitrary index
$\alpha \in {\mathcal A}$.
Let $p \colon X \to R$ and $q \colon Y_{\alpha} \to R$
be unital $\ast$-homomorphisms into a unital $C^{\ast}$-algebra
$R$ such that $p\circ i_{\alpha} = q\circ f_{\alpha}$.
Consider the homomorphisms $\widetilde{p} =
p\circ \varphi_{X} \colon Z \to R$ and
$\widetilde{q} = q\circ\varphi_{Y}^{\alpha}
\colon Z_{\alpha}\bb K \to R$. Note that
\begin{multline*}
\widetilde{p}\circ s_{\alpha} =
p\circ \varphi_{X}\circ s_{\alpha} =
p\circ i_{\alpha}\circ\varphi_{X}^{\alpha} = q\circ f_{\alpha}\circ\varphi_{X}^{\alpha}= q\circ\varphi_{Y}^{\alpha}\circ\pi_{Z_{\alpha}} = \widetilde{q}\circ\pi_{Z_{\alpha}}.
\end{multline*}

\noindent Now let $\alpha \in {\mathcal A}$. Since, by Lemma \ref{L:pushoutex}, the diagram
(the back face of the above cubic diagram)

\[
\begin{CD}
Z @>\pi_{Z}>> Z\bb K\\
@A{s_{\alpha}}AA @AA{s_{\alpha}\bb\operatorname{id}_{K}}A \\
Z_{\alpha} @>\pi_{Z_{\alpha}}>> Z_{\alpha}\bb K
\end{CD}
\]

\noindent is a pushout, it follows that there exists a unique unital
$\ast$-homomorphism $\widetilde{r} \colon Z\bb K \to R$ such that
$\widetilde{p} = \widetilde{r}\circ\pi_{Z}$ and
$\widetilde{q} = \widetilde{r}\circ
(s_{\alpha}\bb \operatorname{id}_{K})$.
Now let $r = \widetilde{r} \circ\phi_{Y} \colon Y \to R$.
We have
\[ r\circ j_{\alpha} = \widetilde{r}\circ\phi_{Y}\circ j_{\alpha}
= \widetilde{r}\circ (s_{\alpha}\bb \operatorname{id}_{K})
\circ \phi_{Y}^{\alpha} = \widetilde{q}\circ\phi_{Y}^{\alpha} =
q\circ \varphi_{Y}^{\alpha}\circ\phi_{Y}^{\alpha} = q \]

\noindent and
\[r\circ f = \widetilde{r}\circ\phi_{Y}\circ f =
\widetilde{r}\circ\pi_{Z}\circ\phi_{X} =
\widetilde{p}\circ\phi_{X} = p \circ\varphi_{X}\circ\phi_{X} = p .\]

\noindent This simply means that the diagram under consideration
has the corresponding universality property. Finally the
uniqueness of $\widetilde{r}$ guarantees that $r$ is the only
unital $\ast$-homomorphism with the just indicated properties.
This shows that our diagram is pushout and completes the proof of part I.

{\em Part II}. Suppose that we are given direct
$C^{\ast}_{\omega}$-systems
${\mathcal S}_{X} = \{ X_{\alpha}, i_{\alpha}^{\beta}, {\mathcal A}\}$,
${\mathcal S}_{Y} = \{ Y_{\alpha}, j_{\alpha}^{\beta},
{\mathcal A}\}$ and a morphism $\{ f_{\alpha} \colon
X_{\alpha} \to Y_{\alpha}; \alpha \in {\mathcal A}\}
\colon {\mathcal S}_{X} \to {\mathcal S}_{Y}$, satisfying
the above indicated properties. 

Let $\alpha \in {\mathcal A}$. By conditions (d), (e) and Lemma
\ref{L:comp}, the composition $j_{\alpha}\circ f_{\alpha}$ is doubly
projective. Since, by condition (d), the inclusion $i_{\alpha}$ is
doubly projective, it follows from Lemmas \ref{L:inj} and \ref{L:out},
that $f$ is also doubly projective. By condition (c) and Lemma
\ref{L:sepsep}, the homomorphism $f_{\alpha}$ has a separable type.
Finally, by condition (f) and Lemma \ref{L:seppushout},
$f$ also has a separable type.
\end{proof}

\begin{cor}\label{C:q}
Let $f \colon X \to Y$ be a doubly projective homomorphism of unital
projective $C^{\ast}$-algebras. If $f$ has a separable type, then
there exists a pushout
\[
\begin{CD}
X @>f>> Y\\
@A{p}AA @AA{q}A\\
X_{0} @>f_{0}>> Y_{0},
\end{CD}
\]

\noindent where $X_{0}$ and $Y_{0}$ are separable unital projective
$C^{\ast}$-algebras and the homomorphisms $i_{0} \colon X _{0} \to Y_{0}$,
$p \colon X_{0} \to X$ and $q \colon Y_{0} \to Y$ are doubly projective.
\end{cor}

\begin{rem}\label{R:full}
Combining methods of proofs of Theorems \ref{T:charact} and
\ref{T:q} it is possible to obtain a characterization of arbitrary
(not necessarily of a separable type) doubly projective
homomorphisms of unital $C^{\ast}$-algebras. This characterization
is recorded in Theorem \ref{T:main}. We only note here that the
sufficiency follows from Proposition \ref{P:doubly}.
\end{rem}

\begin{thm}\label{T:main}
A unital $\ast$-homomorphism $f \colon X \to Y$ of projective
unital $C^{\ast}$-algebras is doubly projective if and only
if there exist well ordered continuous direct systems
${\mathcal S}_{X} = \{ X_{\alpha}, i_{\alpha}^{\alpha +1},
\tau\}$, ${\mathcal S}_{Y} = \{ Y_{\alpha},
j_{\alpha}^{\alpha +1}, \tau\}$ and a morphism
$\{ f_{\alpha}; \tau\} \colon {\mathcal S}_{X}
\to {\mathcal S}_{Y}$
satisfying the following conditions:
\begin{itemize}
\item[(a)]
$X = \varinjlim{\mathcal S}_{X}$,
$Y = \varinjlim{\mathcal S}_{Y}$ and
$f = \varinjlim\{ f_{\alpha} ; \tau\}$.
\item[(b)]
$C^{\ast}$-algebras $X_{0}$ and $Y_{0}$ are separable projective
and the homomorphism
$f_{0} \colon X_{0} \to Y_{0}$ is doubly projective.
\item[(c)]
$C^{\ast}$-algebras $X_{\alpha}$ and $Y_{\alpha}$ are projective
and the homomorphism $f_{\alpha} \colon X_{\alpha}\to Y_{\alpha}$
is doubly projective, $\alpha < \tau$.
\item[(d)]
All short injections $i_{\alpha}^{\alpha +1} \colon
X_{\alpha} \to X_{\alpha +1}$ and $j_{\alpha}^{\alpha +1}
\colon Y_{\alpha} \to Y_{\alpha +1}$ are doubly projective
and have a separable type.
\item[(e)]
All adjacent square diagrams

\[
\begin{CD}
X_{\alpha +1} @>f_{\alpha +1}>> Y_{\alpha +1}\\
@A{i_{\alpha}^{\alpha +1}}AA @AA{j_{\alpha}^{\alpha +1}}A\\
X_{\alpha} @>f_{\alpha}>> Y_{\alpha}
\end{CD}
\]

\noindent are doubly projective and their characteristic
homomorphisms have
separable type.
\item[(f)]
If the homomorphism $f$ itself has a separable type,
then all the square
diagrams indicated in (d) are pushouts.
\end{itemize}
\end{thm}


\end{document}